\documentclass[12pt,a4paper,preprintnumbers,showkeys]{revtex4-1} 

\usepackage{tikz}
\usepackage[fleqn]{amsmath} 
\usepackage{amssymb} 
\usepackage{mathrsfs} 
\usepackage{bm}
\usepackage{graphicx} 
\usepackage{subfigure}  
\usepackage{pict2e} 
\usepackage{ulem} 
\usepackage{cases}
\usepackage[pdfstartview=FitH,
            bookmarksnumbered=true,
            bookmarksopen=true,
            colorlinks=true, 
            pdfborder=001,   
            linkcolor=blue,  
            citecolor=blue,   
            urlcolor=blue    
            ]{hyperref}      
\usepackage[left=2.5cm, right=2.5cm, top=2.5cm, bottom=2.2cm]{geometry}
\newtheorem{definition}{Definition}[section]
\newtheorem{lemma}{Lemma}[section]
\newtheorem{theorem}{Theorem}[section]

\begin{document}

\title{Constructing the $r$-uniform supertrees with \\ the same spectral radius and matching energy}

\author{Jin-Xiu Zhou, Wen-Huan Wang}\thanks{Corresponding author. Email: whwang@shu.edu.cn}

\affiliation{Department of Mathematics, Shanghai University, Shanghai 200444, China}

\date{\today}

\begin{abstract}

  An $r$-uniform supertree is a connected and acyclic hypergraph of which each edge has $r$ vertices, where $r\geq 3$.
   \textcolor{black}{We propose the concept of matching energy for an $r$-uniform hypergraph, which is  defined as  the sum of the absolute value of all the eigenvalues of its matching polynomial}.
    With the aid of the matching polynomial of an $r$-uniform supertree,
   three pairs of $r$-uniform supertrees with the same spectral radius \textcolor{black}{and the same matching energy} are  constructed, and two  infinite families of $r$-uniform supertrees with the same spectral radius \textcolor{black}{and the same  matching energy} are characterized.
    Some known results about the graphs with the same spectra regarding to their adjacency matrices can be  naturally deduced from our new results.
      \end{abstract}

\keywords{Supertree; Matching polynomial; Spectral radius; \textcolor{black}{Matching energy}}


\maketitle


\section{Introduction}

  Let   $\mathbb{C}$  and $\mathbb{R}$ be the sets of  complex and real numbers, respectively. Let $r$ and $s$ be two positive integers not less than 2 and $[s] = \{1,\ldots,s\}$.
   We denote by $\mathcal{A}=(a_{i_1i_2\cdots i_r})$ a  real tensor (or hypermatrix)  of order $r$ and dimension $s$, which  is a multi-dimensional array
   with entries $a_{i_1i_2\cdots i_r}\in \mathbb{R}$, where  $i_1,i_2,\cdots, i_r\in [s]$ and $r\geq 2$.
   If $r=2$, then $\mathcal{A}$ is a matrix.
    If $a_{i_1i_2\cdots i_r}=1$ when  $i_{1}=i_{2}=\cdots=i_{r}$ and $a_{i_1i_2\cdots i_r}=0$ otherwise, then  $\mathcal{A}$ is the identity  tensor.
        Let   $ \bm{x}=(x_1,x_2,\ldots,x_s)^{\textrm{T}}\in \mathbb{C}^s$
  be an  $s$-dimensional complex column vector    and 
  $ \bm{x}^{[r]}=(x^{r}_1,x^{r}_2,\cdots,x^{r}_s)^{\textrm{T}}$.
  Then
  $\mathcal{A}\bm{x}$
  is a vector in $\mathbb{C}^s$ whose $i$-th component is given by
  \begin{align}\label{444444}
  (\mathcal{A}\bm{x})_i=\sum^s_{i_2,\ldots,i_r=1}a_{ii_2\cdots i_r}x_{i_2}\cdots x_{i_r},~\mbox{for~each}~i\in [s].
  \end{align}
 Furthermore, we have
     \begin{equation}\label{2.2}
  \bm{x}^{\textrm{T}} (\mathcal{A}\bm{x})=\sum_{i_1,i_{2},\ldots,i_{r}=1}^{s}a_{i_1 i_{2}\ldots i_{r}}x_{i_{1}}\cdots x_{i_{r}}.
  \end{equation}




  The concept of tensor eigenvalues and the spectra of tensors was introduced by   Qi \cite{Qi2005}  and Lim \cite{Lim2005}  in 2005 independently as follows.
     If there exist  $\lambda\in \mathbb{C}$ and  $\bm{x}\in \mathbb{C}^s$ satisfying $\mathcal{A}\bm{x}=\lambda \bm{x}^{[r-1]}$, namely,  $(\mathcal{A}\bm{x})_{i}=\lambda x_i^{r-1}$ for  any $i\in[s]$,  then   $\lambda$ is called an eigenvalue of  $\mathcal{A}$ and $\bm{x}$ an
 eigenvector of  $\mathcal{A}$ corresponding to  $\lambda$.
  The resultant of the $s $ element homogeneous equations $\mathcal{A} \bm{x}=0$ is called the determinant of  $\mathcal{A}$ and is denoted by $\operatorname{det}(\mathcal{A})$.
  The characteristic polynomial of $\mathcal{A}$ is defined as $\Phi(\mathcal{A},x)=\operatorname{det}\left(x \mathcal{I}-\mathcal{A}\right)$, where  $\mathcal{I}$ is the  identity  tensor of order $r$ and dimension $s$. The eigenvalues of $\mathcal{A}$ are the roots of $\Phi(\mathcal{A},x)$ \cite{Shao2013}.
 The (multi)-set of all the roots of $\Phi(\mathcal{A},x)$ (counting multiplicities), denoted by $\operatorname{Spec}\mathcal{A}$, is called the spectra of $\mathcal{A}$.

  A hypergraph $\mathcal{H}$ is an ordered pair $(V(\mathcal{H}), E(\mathcal{H}))$, where $V(\mathcal{H})=[s]$ is the set of vertices of $\mathcal{H}$   and
 $E(\mathcal{H}) \subseteq P([s])$ the set of edges of $\mathcal{H}$ with $P([s])$ being the power set of $[s]$.
 If each edge $e$ of $E(\mathcal{H})$ has $r$ vertices ($r \geq 2$),
 then  $\mathcal{H}$ is an $r$-uniform hypergraph.
 If $r=2$, then $\mathcal{H}$ reduces to an ordinary graph and \textcolor{black}{we denote it by $H$}.
   A hypergraph  $\mathcal{H}$ is called linear  if
   any two edges of $\mathcal{H}$  intersect on at most one common vertex.
  If $\mathcal{H}$  does not contain cycles, then $\mathcal{H}$ is  acyclic or a superforest.
  If $\mathcal{H}$  is connected and acyclic, then $\mathcal{H}$ is a supertree   \cite{Li2016d}.
   In this paper, we consider  $r$-uniform supertrees.

 Let $\mathcal{H}=(V(\mathcal{H}), E(\mathcal{H}))$ be an $r$-uniform hypergraph on $s$ vertices.
  The adjacency tensor of $\mathcal{H}$ is the $r$-ordered  and $s$-dimensional tensor  $\mathcal{A}(\mathcal{H})=\left(a_{i_{1} i_{2} \cdots i_{r}}\right)$,
  where $a_{i_{1} i_{2} \cdots i_{r}}= \frac{1}{(r-1) !}$ if $\{i_{1}, i_{2}, \ldots, i_{r}\} \in E(\mathcal{H})$ and
  0 otherwise \cite{Cooper2012}.
 The spectral radius of $\mathcal{H}$, denoted by $\rho(\mathcal{H})$,
  is defined as the maximum modulus of all the eigenvalues of the characteristic polynomial $\Phi(\mathcal{A}(\mathcal{H}),x)$.

   Let $\bm{M}$ be a matrix.
   It may stand for the adjacency matrix, the Laplacian matrix, the signless
   Laplacian matrix, and the distance matrix, etc.
  Two graphs are said to be $\bm{M}$-cospectral if they have the same $\bm{M}$-spectra, where $\bm{M}$-spectra of a graph is
  the (multi)-set of all the eigenvalues of its corresponding $\bm{M}$ matrix.
   \textcolor{black}{Similarly, two hypergraphs are said to be adjacency cospectral, if their adjacency tensors have the same characteristic polynomial. A graph $G$ (a hypergraph $\mathcal{G}$, respectively) is determined by its $\bm{M}$-spectra (spectra, respectively) if there does not exist other non-isomorphic graph $H$ (hypergraph $\mathcal{H}$, respectively)
    such that $H$ and $G$ ($\mathcal{G}$ and $\mathcal{H}$, respectively) are $\bm{M}$-cospectral (cospectral, respectively).
  }

   Which graphs are determined by their spectra?
   G\"{u}nthard and Primas \cite{gunthard1956} posed this fundamental problem in 1956 in the context of H\"{u}ckel's theory in chemistry.
   Constructions of cospectral non-isomorphic graphs have implications on the complexity of the graph isomorphism
    problem and reveal which graph properties cannot be deduced from the spectra of graphs.
    Therefore, it  can  help researchers understand the above  question.
    The construction of cospectral graphs attracted many researcher's attention and has been studied extensively.

 When $\bm{M}$ is an adjacency matrix, many results about \textcolor{black}{the $\bm{M}$-cospectral graphs}
  have been obtained.
 In the 1960s,  Van Lint and Seidel \cite{Lint1966} introduced the Seidel switching for
 constructing families of cospectral graphs. 
  By using the Seidel switching,
 recently Seress \cite{Seress2000}  constructed an infinite
  family of cospectral eight regular graphs.
   Godsil and McKay \cite{godsil1982} further developed this
 concept and introduced the Godsil--McKay switching.
  Bl\'{a}zsik et al.\  \cite{Blazsik2015} used the Godsil--McKay  switching to
 construct two cospectral regular graphs such that one has a perfect matching while the  other does not have any perfect matching.
 Langberg and Vilenchik \cite{Langberg2018}    presented a new method which was based on bipartite graph product to construct an infinite family of cospectral graphs.
   Qiu et al.\ \cite{Qiu2020} 
    constructed an infinite
  family of cospectral graphs by using  new methods. 
   For oriented graphs and signed graphs, Belardo et al.\  \cite{2021-Belardo-p2717} extended the Godsil--McKay switching to signed graphs,  and
    built pairs of cospectral switching nonisomorphic signed graphs and
  Stani\'{c} \cite{Stanic2020} obtained infinite families of cospectral regular signed graphs and cospectral bipartite regular oriented graphs.

  When $\bm{M}$ is a  Laplacian  matrix and a distance matrix, for the construction of $\bm{M}$-cospectral  graphs,
  the readers can refer to Refs. 
  \cite{Heysse2017,wen2018,Ajmal2021}.

  To the  authors' best knowledge, the result on the construction of $E$-cospectral  hypergraphs is as follows.
  Let $\mathcal{A}$ be a tensor of order $r\geq 2$ and  dimension $s\geq 2$.
    If there exists a nonzero vector $\bm{x}\in \mathbb{C}^s$ such that  $\mathcal{A}\bm{x}=\lambda \bm{x}$ and $\bm{x}^{\top}\bm{x}=1$, then
     $\lambda\in \mathbb{C}^s$ is called an $E$-eigenvalue of $\mathcal{A}$.
  The $E$-eigenvalues of $\mathcal{A}$ are the roots of the $E$-characteristic polynomial $\phi_{\mathcal{A}}(\lambda)$ of $\mathcal{A}$ (see \cite{2007-Qi-p1363} for the definition of $\phi_{\mathcal{A}}(\lambda)$).
  Recently, Bu et al.\ \cite{Bu2014} deduced a method of constructing $E$-cospectral hypergraphs
   and obtained some hypergraphs which are determined by their spectra.
   However,  as we  know, the calculation of the characteristic polynomial  of  the adjacency tensor of \textcolor{black}{hypergraph} is NP-hard in any field \cite{Grenet2010}.
   Therefore, it is difficult for us to use the characteristic polynomial
   of \textcolor{black}{hypergraph} to study the cospectral hypergraph.
  Since the spectral radius of \textcolor{black}{hypergraph} is of practical significance {\cite{2017-ChenOuyang-p141},
  the characterization of the $r$-uniform \textcolor{black}{hypergraph} with the extremal spectral radius is interesting, and
  a lot of  results have been obtained.
   The interested readers can refer to Refs.\ \cite{Li2016d,2020-Wang-p,2019-Xiao-p1392}.

    Let $\mathcal{H}$ be  an $r$-uniform hypergraph. 
   The number of $k$-matchings in $\mathcal{H}$, denoted by  $m(\mathcal{H},k)$, is the number of selections of $k$ independent edges in $\mathcal{H}$, where $k\geq 0$.
   For the sake of consistence, let $m(\mathcal{H},0)=1$.
   The matching number $\nu(\mathcal{H})$ of $\mathcal{H}$ is the maximum cardinality of a matching in $\mathcal{H}$.
   The matching polynomial of $\mathcal{H}$, denoted by
   $ \varphi(\mathcal{H}, x)=\sum_{k=0}^{\nu(\mathcal{H})}(-1)^{k} m(\mathcal{H}, k) x^{(\nu(\mathcal{H})-k) r}$, was first introduced by Zhang et al.\ \cite{Zhang2017a} when they studied  the spectra of $r$-uniform supertrees.
   In order to  guarantee that the matching polynomials of the $r$-uniform hypergraphs with $n$ vertices have the same degree, $ \varphi(\mathcal{H}, x)$ is  redefined by  Su et al.\ \cite{2018-Su-p4} as $ \varphi(\mathcal{H}, x)=\sum_{k \geq 0}(-1)^{k} m(\mathcal{H}, k) x^{n-k r}$.
   \textcolor{black}{It is noted that when $r=2$, $\mathcal{H}$ is an ordinary graph (denoted by $H$) and $ \varphi(H, x)$ is the matching polynomial of $H$.}

    \textcolor{black}{The matching energy of an ordinary graph $H$, denoted by $ME(H)$,
    was proposed by  Gutman and Wagner \cite{Gutman-2012-p2177}  and it was defined as the sum of the absolute value of all the eigenvalues of $ \varphi(H, x)$.
    Gutman and Wagner \cite{Gutman-2012-p2177} pointed out that $ME(H)$ is a quantity which has a close relationship with chemical applications and it can be traced back to the 1970s.
     For more details about the matching energy, one can refer to \cite{Gutman-2012-p2177}.
   In this paper,  we  extend the definition of the matching energy of a graph $H$ to an $r$-uniform hypergraph $\mathcal{H}$ as follows.
   Similarly, we define the matching energy of $\mathcal{H}$ as  the sum of the absolute value of all the eigenvalues of $ \varphi(\mathcal{H}, x)$,  and denote by $ME(\mathcal{H})$ the matching energy of $\mathcal{H}$.
   We expect that $ME(\mathcal{H})$ can be applied in chemistry as $ME(H)$ does.}

  Motivated by the above-mentioned results, in this paper,  we will study the $r$-uniform supertrees with the same spectral radius of their adjacency tensors and \textcolor{black}{the same matching energy}.  \textcolor{black}{Hereinafter, for simplicity, the $r$-uniform supertrees with the same spectral radius and the same matching energy is abbreviated to the $r$-uniform supertrees with the same SR and   ME, where SR and  ME stand for the spectral radius and the matching energy, respectively. }
  The mail tool used is the matching polynomial of the $r$-uniform supertrees.
  This paper is organized as follows.
  In Section 2, some basic definitions and necessary lemmas are introduced.
  In Sections 3 and 4, the first and the second pairs of  $r$-uniform supertrees with the same \textcolor{black}{SR and ME} are characterized
  (as shown in Theorems \ref{c2w} and \ref{c4}, respectively) and  two infinite families of  $r$-uniform supertrees with the same \textcolor{black}{SR and ME} are constructed (as shown in Theorems \ref{wuxian} and \ref{c4}).
   Three pairs of graphs which are
    $\bm{M}$-cospectral are deduced (as shown in Theorems \ref{c1}, \ref{c11} and \ref{c3s}) in Sections 3 and 4, where $\bm{M}$ is the adjacency matrix.
  In Section 5, we characterize the third  pair of $r$-uniform supertrees with the same \textcolor{black}{SR and ME}  (as  shown in
  Theorem \ref{c5})
  and get a graph which is not determined by its  spectra of its adjacency matrix (as  shown in
  Theorem \ref{c6}).

\section{Preliminary}

 In this section,  some notations and  necessary lemmas are introduced.

 Let $\mathcal{H}$ be a hypergraph and  $v$ be a vertex of $\mathcal{H}$.
 Let $E_{\mathcal{H}}(v)$ be the set of  the edges incident with $v$ and
  $d_{\mathcal{H}}(v)$ the degree  of $v$.
 Namely, $d_{\mathcal{H}}(v)=|E_{\mathcal{H}}(v)|$.
  For $e=\{u_1,\ldots,u_r\}\in E(\mathcal{H})$, if $d_{\mathcal{H}}(u_1)\geq2$ and  $d_{\mathcal{H}}(u_i)=1$ for $2\leq i\leq r$, then  we say that $e$ is a  pendent edge at $u_1$ of $\mathcal{H}$.  If $d_{\mathcal{H}}(v)=1$ and  $v$ is incident with a pendent edge of $\mathcal{H}$, then $v$ is said to be a pendent vertex.

  Let $\mathcal{H}-v$ be the hypergraph obtained from $\mathcal{H}$ by deleting $v$ together with all the edges in $E_{\mathcal{H}}(v)$.
  For $e\in E(\mathcal{H})$ and $V(e)=\{u_1,\cdots,u_r\}$,
     $\mathcal{H}-V(e)$ is the hypergraph obtained from $\mathcal{H}$ by deleting all the vertices in $V(e)$.
    For a subset $E'\subseteq E(\mathcal{H}) $  in $\mathcal{H}$,  $\mathcal{H}\backslash E'$ is the hypergraph obtained from $\mathcal{H}$ by deleting all the edges in $E'$. Namely, $\mathcal{H}\backslash E'=(V(\mathcal{H}), E(\mathcal{H})\setminus E')$.
   If $E'=\{e\}$, then we write $\mathcal{H}\backslash E'$ as $\mathcal{H}\backslash e$.
    Let $N_{k}$ be the set of $k$ isolated vertices, where $k \geq 1$.
  Let  $\mathcal{G}\cup\mathcal{H}$ be the  union of $\mathcal{G}$ and $\mathcal{H}$, where $\mathcal{G}$ and $\mathcal{H}$ are two disjoint hypergraphs.
   If  $V'\subseteq V(\mathcal{H})$ and $E'\subseteq E(\mathcal{H})$, then
  $\mathcal{H'}=(V',E')$ is a partial hypergraph of $\mathcal{H}$.
 Furthermore, if $\mathcal{H'}\neq\mathcal{H}$, then
  $\mathcal{H'}$ is a proper partial hypergraph of $\mathcal{H}$.

   Let $G$ be a graph. We denote by   $\phi(G,x)$ and  $\varphi(G,x)$  the characteristic polynomial and the matching polynomial of   $G$, respectively.

\begin{lemma}\label{tp}(\cite{godsil1993})
  If $G$ is a forest, then $\phi(G,x)=\varphi(G,x)$.
  \end{lemma}

}

  Friedland et al.\ \cite{Friedland2013} defined the nonnegative weakly irreducible tensor
  and Yang et al.\ \cite{-N-p} restated it as follows.

\begin{definition}\label{definition2.3} \cite{-N-p}
 Let $\mathcal{A}=(a_{i_1i_2\cdots i_r})$ be a nonnegative tensor of order $r$ and dimension  $s$.
 If for any nonempty proper index subset $I \subset [s]$, there is at least an entry $a_{i_1i_2\cdots i_r}>0$,
 where $i_{1}\in I$ and at least an $i_{j}\in {{[s]\setminus I}}$ for $j=2,3, \ldots, r$,
 then $\mathcal{A}$ is called a nonnegative weakly irreducible tensor.
  \end{definition}

     It was proved that an $r$-uniform hypergraph $\mathcal{H}$ is connected if and only if
        its adjacency tensor $\mathcal{A}(\mathcal H)$ is weakly irreducible (see \cite{Friedland2013} and \cite{-N-p}).

 \begin{lemma}\label{pf}\cite{Friedland2013,Yang2010}
 Let $\mathcal{A}$ be a nonnegative tensor of order $r$ and dimension $s$, where $r\geq 2$.
 Then we have the following statements.

 (i). $\rho(\mathcal{A})$ is an eigenvalue of $\mathcal{A}$ with a nonnegative eigenvector $\bm{x}\in \mathbb{R}^{s}_{+}=\{x\in \mathbb{R}^{s}\mid x\geq 0\}$ corresponding to it.

 (ii). If $\mathcal{A}$ is weakly irreducible,
 then $\rho(\mathcal{A})$ is the only eigenvalue of $\mathcal{A}$ with a positive eigenvector $\bm{x}\in \mathbb{R}^{s}_{++}=\{x\in \mathbb{R}^{s}\mid x> 0\}$, up to a positive scaling coefficient.
 \end{lemma}

\begin{lemma}\label{sub} (\cite{Cooper2012,fan2015})
 Suppose that $\mathcal{H}$ is a uniform hypergraph, and $\mathcal{H}^{\prime}$ is a partial hypergraph of $\mathcal{H} .$ Then $\rho\left(\mathcal{H}^{\prime}\right) \leq \rho(\mathcal{H})$.  Furthermore, if $\mathcal{H}$ is connected and $\mathcal{H}^{\prime}$ is a proper partial hypergraph, we have $\rho\left(\mathcal{H}^{\prime}\right)<\rho(\mathcal{H})$.
 \end{lemma}

 \begin{lemma} (\cite{Cooper2012})
  Let $\mathcal H$ be an $r$-uniform hypergraph that is the disjoint union of hypergraphs $\mathcal{H}_{1}$ and $\mathcal{H}_{2}$. Then as sets, $\operatorname{Spec}(\mathcal H)=\operatorname{Spec}\left(\mathcal{H}_{1}\right) \cup \operatorname{Spec}\left(\mathcal{H}_{2}\right)$. Considered as multisets, an eigenvalue $\lambda$ with multiplicity $m$ in $\operatorname{Spec}\left(\mathcal{H}_{1}\right)$ contributes $\lambda$ to $\operatorname{Spec}(\mathcal H)$ with multiplicity $m(r-1)^{\left|H_{2}\right|}$.
  \label{bing}
 \end{lemma}

  A totally nonzero eigenvalue of hypergraph $\mathcal{H}$ is a nonzero eigenvalue and all the entries of the eigenvectors corresponding to it are nonzero.


 \begin{lemma}\label{redf}(\cite{Zhang2017a,2018-Su-p4})
 $\lambda$ is a totally nonzero eigenvalue of an $r$-uniform supertree $\mathcal{H}$ with $n\geq 3$ vertices if and only if it is a root of the matching polynomial
\begin{align*}
 \varphi(\mathcal{H}, x)=\sum_{k \geq 0}(-1)^{k} m(\mathcal{H}, k) x^{n-k r} .
\end{align*}
  \end{lemma}

 \begin{lemma} \label{su} (\cite{2018-Su-p4})
 Let $\mathcal{G}$ and $\mathcal{H}$ be two $r$-uniform hypergraphs. Then the following statements hold.\\
 (a) $\varphi(\mathcal{G} \cup \mathcal{H}, x)=\varphi(\mathcal{G}, x) \varphi(\mathcal{H}, x)$.\\
 (b) If $u \in V(\mathcal{G})$ and $I=\left\{i \mid e_{i} \in E_{\mathcal{G}}(u)\right\}$, then for any $J \subseteq I$, we have
\begin{align*}
 \varphi(\mathcal{G}, x)=\varphi\left(\mathcal{G} \backslash\left\{e_{i}: i \in J\right\}, x\right)-\sum_{i \in J} \varphi\left(\mathcal{G}-V\left(e_{i}\right), x\right),
 \end{align*}
\begin{align*}
 \varphi(\mathcal{G}, x)=x \varphi(\mathcal{G}-u, x)-\sum_{e \in E_{\mathcal{G}}(u)} \varphi(\mathcal{G}-V(e), x) .
\end{align*}
  \end{lemma}
 \begin{lemma} \label{large}
 Let $\mathcal{H}$ be an $r$-uniform hypergraph with $n\geq 3$ vertices.
  Then $\rho(\mathcal{H})$  is the largest root of $\varphi(\mathcal{H}, x)=\sum\limits_{k \geq 0}(-1)^{k} m(\mathcal{H}, k) x^{n-k r}$.
  \end{lemma}
 \noindent\textbf{Proof}. By  Lemma \ref{pf}, $\rho(\mathcal{H})$ is a totally nonzero eigenvalue of $\mathcal{H}$.
 The set of totally nonzero eigenvalues of $\mathcal{H}$ is denoted by $M=\{\lambda_1,\lambda_2,\cdots ,\lambda_l,\rho(\mathcal{H})\}$, where $l$ is a positive integer.
  Without loss of generally,  we suppose $|\lambda_1|\leq|\lambda_2|\leq\cdots\leq |\lambda_l|\leq\rho(\mathcal{H})$.
 Let $N=\{\mu_1,\mu_2,\cdots ,\mu_{l'}\}$ be the set of the nonzero roots of $\varphi(\mathcal{H}, x)$, where
   $|\mu_1|\leq|\mu_2|\leq\cdots \leq|\mu_{l'}|$ and $l'$ is a positive integer.
 It follows from Lemma \ref{redf} that   $\rho(\mathcal{H})=\mu_{l'}$.
 ~~$\Box$
%
%
%
%
%
%

 \section{The first pair of  $r$-uniform supertrees with the same spectral radius and matching energy}

 In this section, we first characterize the first  pair of $r$-uniform supertrees with the same SR \textcolor{black}{and ME} in Theorem \ref{c2w}.
 Then, from Theorem \ref{c2w}, we obtain an infinite family of  $r$-uniform supertrees with the same SR \textcolor{black}{and  ME} in Theorem \ref{wuxian}.
 Furthermore,   
 an example  in Theorem \ref{wth4} is given to show how to use  Theorem \ref{c2w} to determine whether two $r$-uniform supertrees have the same SR \textcolor{black}{and ME} or not.
 Finally, from  our new results, a known pair of graphs and a known infinite family  of graphs which are  $\bm{M}$-cospectral are naturally deduced (as shown in Theorems \ref{c1} and  \ref{c11}, respectively), where $\bm{M}$ is the adjacency matrix.


  Let $\mathcal{G}$ and $\mathcal{H}$ be two $r$-uniform supertrees whose vertex sets are disjoint with $u\in V(\mathcal{G})$ and $v\in V(\mathcal{H})$, where $r \geq 2$. 
   We denote by $\mathcal{G}(u,v)\mathcal{H}$ the supertree obtained from $\mathcal{G}$ and $\mathcal{H}$ by identifying $u$    with  $v$.
  To obtain our results, we introduce Lemma \ref{th1} as follows.

 \begin{lemma}\label{th1}
  Let $\mathcal{G}$, $\mathcal{H}$, and $\Gamma$ be three $r$-uniform supertrees, where   $\mathcal{G}$ and  $\mathcal{H}$ have the same number of vertices and $r\geq 2$.
  Let $u\in V(\mathcal{G})$ and   $v\in V(\mathcal{H})$.
  If $\varphi(\mathcal{G},x)=\varphi(\mathcal{H},x)$ and  $\varphi(\mathcal{G}-u,x)=\varphi(\mathcal{H}-v,x)$, then  for any  $w\in V(\Gamma)$, we have $\varphi(\mathcal{G}(u,w)\Gamma,x)=\varphi(\mathcal{H}(v,w)\Gamma,x)$.
  \end{lemma}

 \noindent\textbf{Proof}.
 In $\mathcal{G}(u,w)\Gamma$, let $q$ be the vertex  $u$ of $G$ (namely $w$ of $\Gamma$).
 For  simplicity, let $\mathcal{G}(u,w)\Gamma=\mathcal{G}\cdot\Gamma $ and $\mathcal{H}(v,w)\Gamma=\mathcal{H}\cdot\Gamma$.
 By  Lemma \ref{su}(b), we get
 \begin{align}\label{11}
  \varphi(\mathcal{G}\cdot\Gamma,x)=x\varphi(\mathcal{G}\cdot\Gamma-q,x)-\sum_{e\in E_{\mathcal{G}\cdot\Gamma}(q)}\varphi(\mathcal{G}\cdot\Gamma-V(e),x).
 \end{align}
 Since $\mathcal{G}\cdot\Gamma-q\cong \left(\mathcal{G}-u\right)\cup\left(\Gamma-w\right)$
 and $E_{\mathcal{G}\cdot\Gamma}(q)=E_{\mathcal{G}}(u)\cup E_{\Gamma}(w)$,  by Lemma \ref{su}(a), we get
\begin{align}\label{22}
   &\varphi(\mathcal{G}\cdot\Gamma,x)=x\varphi(\mathcal{G}-u,x)\varphi(\Gamma-w,x)\nonumber\\
    &- \sum\limits_{e\in E_{\mathcal{G}}(u)}\varphi\left(\mathcal{G}-V(e),x\right)\varphi(\Gamma-w,x)-\sum\limits_{e\in E_{\Gamma}(w)}\varphi\left(\Gamma-V(e),x\right)\varphi(\mathcal{G}-u,x).
\end{align}
 Furthermore, by Lemma \ref{su}(b), we obtain 
 \begin{align}\label{44}
  &\varphi(\mathcal{G}\cdot\Gamma,x)
    =x\varphi(\mathcal{G}-u,x)\varphi(\Gamma-w,x)\nonumber\\
    &+\left[\varphi(\mathcal{G},x)-x\varphi(\mathcal{G}-u,x)\right]\varphi(\Gamma-w,x)
  +\left[\varphi(\Gamma,x)-x\varphi(\Gamma-w,x)\right]\varphi(\mathcal{G}-u,x).
 \end{align}
 Therefore, by simplification, we get
 \begin{align}
  \varphi(\mathcal{G}\cdot\Gamma,x)=\varphi(\mathcal{G},x)\varphi(\Gamma-w,x)+
   \varphi(\mathcal{G}-u,x)\varphi(\Gamma,x)
    -x\varphi(\mathcal{G}-u,x)\varphi(\Gamma-w,x).\label{55}
\end{align}
 Similarly, we get
\begin{equation}\label{66}
  \varphi(\mathcal{H}\cdot\Gamma,x)
  =\varphi(\mathcal{H},x)\varphi(\Gamma-w,x)+\varphi(\mathcal{H}-v,x)\varphi(\Gamma,x)
    -x\varphi(\mathcal{H}-v,x)\varphi(\Gamma-w,x).
\end{equation}
 If $\varphi(\mathcal{G},x)=\varphi(\mathcal{H},x)$ and $\varphi(\mathcal{G}-u,x)=\varphi(\mathcal{H}-v,x)$ hold, then by comparing \eqref{55} and \eqref{66}, we get
 $\varphi(\mathcal{G}\cdot\Gamma,x)=\varphi(\mathcal{H}\cdot\Gamma,x)$.
 ~~$\Box$

  By  Lemmas  \ref{large} and \ref{th1},  we can directly get Theorem \ref{c2w}.

 \begin{theorem} \label{c2w}
  Let $\mathcal{G}$, $\mathcal{H}$ and $\Gamma$ be three $r$-uniform  supertrees,
  where  $\mathcal{G}$ and  $\mathcal{H}$ have the same number of vertices and $r\geq 3$.
    Let $u\in V(\mathcal{G})$ and $v\in V(\mathcal{H})$.
  If $\varphi(\mathcal{G},x)=\varphi(\mathcal{H},x)$ and $\varphi(\mathcal{G}-u,x)=\varphi(\mathcal{H}-v,x)$, then
  for any $w\in V(\Gamma)$, we have $\rho(\mathcal{G}(u,w)\Gamma)=\rho(\mathcal{H}(v,w)\Gamma)$ and \textcolor{black}{$ME(\mathcal{G}(u,w)\Gamma)=ME(\mathcal{H}(v,w)\Gamma)$}.
  \end{theorem}


  Let $\mathcal{G}$, $\mathcal{H}$ and $\Gamma$ be three $r$-uniform  supertrees, where $r\geq 2$.
  Let  $u\in V(\mathcal{G})$,  $v\in V(\mathcal{H})$, and $w\in V(\Gamma)$.
  Let $m$, $n$, $a$, and $b$ be four positive integers.
    For simplicity, we denote $\underbrace{G\cup \cdots \cup G}_{m}$ by $m\mathcal{G}$.
   Let $\mathcal{G}_u^m$ be the hypergraph obtained from $m\mathcal{G}$ by coalescing $u$ such that the $m$ copies of  $\mathcal{G}$ share a common vertex $u$.
  Similarly, $\mathcal{H}_v^n$ is  defined as that of $\mathcal{G}_u^m$.
  We denote by $\mathcal{G}_u^m \cdot \mathcal{H}_v^n$ the hypergraph obtained from $\mathcal{G}_u^m$  and  $\mathcal{H}_v^n$  by identifying $u$ of $\mathcal{G}_u^m$   with $v$ of $\mathcal{H}_v^n$. In particular, $\mathcal{G}_u^1\cong \mathcal{G}_u$.
  Let
  $\mathcal{G}_u^{a+b}=\mathcal{G}_u^a\cdot\mathcal{G}_u^b$ and $\mathcal{H}_v^{a+b}=\mathcal{H}_v^a\cdot\mathcal{H}_v^b$.
  The hypergraph $\mathcal{G}(u,v)\mathcal{H}(v,w)\Gamma$  is  obtained from $\mathcal{G}$, $\mathcal{H}$ and $\Gamma$ by identifying $u$, $v$ and $w$.
  If $r=2$, we write  $\mathcal{G}_u^m$ and $\mathcal{G}_u^m \cdot \mathcal{H}_v^n$ as   $G_u^m$ and $G_u^m \cdot H_v^n$, respectively.
  Obviously, both of them are graphs.  In particular, $G_u^1\cong G_u$.
 Let   $G_u^{a+b}=G_u^a\cdot G_u^b$ and $H_v^{a+b}=H_v^a\cdot H_v^b$.


 \begin{theorem} \label{wuxian}
 Let $\mathcal{G}$ and $\mathcal{H}$ be two $r$-uniform  supertrees with the same number of vertices, where $r\geq 3$.
 Let  $u\in V(\mathcal{G})$,  $v\in V(\mathcal{H})$ and  $m$ be a positive  integer.
 If $\varphi(\mathcal{G},x)=\varphi(\mathcal{H},x)$ and $\varphi(\mathcal{G}-u,x)=\varphi(\mathcal{H}-v,x)$,
 then we have
 \textcolor{black}{(i). $\rho(\mathcal{H}_v^m)=\rho(\mathcal{H}_v^{m-1}\cdot\mathcal{G}_u)
 =\rho(\mathcal{H}^{m-2}_v\cdot\mathcal{G}_u^2)=\cdots
 =\rho(\mathcal{H}_v\cdot\mathcal{G}_u^{m-1})
 =\rho(\mathcal{G}_u^m)$; (ii). $ME(\mathcal{H}_v^m)=ME(\mathcal{H}_v^{m-1}\cdot\mathcal{G}_u)
 =ME(\mathcal{H}^{m-2}_v\cdot\mathcal{G}_u^2)=\cdots
 =ME(\mathcal{H}_v\cdot\mathcal{G}_u^{m-1})
 =ME(\mathcal{G}_u^m)$}.
\end{theorem}

\noindent\textbf{Proof}.
 Since $\mathcal{H}_v^m=\mathcal{H}_v^{m-1}\cdot \mathcal{H}_v$, if $\varphi(\mathcal{G},x)=\varphi(\mathcal{H},x)$ and $\varphi(\mathcal{G}-u,x)=\varphi(\mathcal{H}-v,x)$,
 then by Lemma \ref{th1}, we have $\varphi(\mathcal{H}_v^m)=\varphi(\mathcal{H}_v^{m-1}\cdot\mathcal{G}_u)$.
 Furthermore, by Lemma  \ref{large}, 
\textcolor{black}{ we have $\rho(\mathcal{H}_v^m)=\rho(\mathcal{H}_v^{m-1}\cdot\mathcal{G}_u)$.
 Similarly, we obtain $\varphi(\mathcal{H}_v\cdot\mathcal{G}_u^{m-1})
 =\varphi(\mathcal{G}_u^m)$ and $\rho(\mathcal{H}_v\cdot\mathcal{G}_u^{m-1})
 =\rho(\mathcal{G}_u^m)$.}
 Next, we only need to prove $\varphi(\mathcal{H}_v^{m-k}\cdot \mathcal{G}_u^{k},x)=\varphi(\mathcal{H}_v^{m-k-1}\cdot \mathcal{G}_u^{k+1},x)$, where $k=1,\cdots,m-2$.

  Let $k=1,\cdots,m-2$. Let $\Gamma=\mathcal{H}_v^{m-k-1}\cdot \mathcal{G}_u^{k}$ and $w$ of $\Gamma$ be $u$ of $\mathcal{G}^{k}$ (namely $v$ of $\mathcal{H}^{m-k-1}$).
  Obviously, $\mathcal{H}_v^{m-k}\cdot \mathcal{G}_u^{k}=H_v \cdot \Gamma_w$  and
     $\mathcal{H}_v^{m-k-1}\cdot \mathcal{G}_u^{k+1}=G_u \cdot \Gamma_w$.
 Since $\varphi(\mathcal{G},x)=\varphi(\mathcal{H},x)$ and $\varphi(\mathcal{G}-u,x)=\varphi(\mathcal{H}-v,x)$, by Lemma \ref{th1}, we obtain
  $\varphi(\mathcal{H}_v^{m-k}\cdot \mathcal{G}_u^{k},x)=\varphi(\mathcal{H}_v^{m-k-1}\cdot \mathcal{G}_u^{k+1},x)$. Furthermore,  by Lemma \ref{large}, we get
   $\rho(\mathcal{H}_v^{m-k}\cdot \mathcal{G}_u^{k},x)=\rho(\mathcal{H}_v^{m-k-1}\cdot \mathcal{G}_u^{k+1},x)$.
  Therefore, we get Theorem \ref{wuxian}(i). \textcolor{black}{By the definition of the matching energy of an $r$-uniform  hypergraph, we obtain
  Theorem \ref{wuxian}(ii).}
   ~~$\Box$


\begin{figure}
\thicklines
\begin{tikzpicture}
\draw[thick,black](2.6,-1.2) ellipse(0.5 and 0.1); 
\draw[thick,black](3.3,-1.2) ellipse(0.5 and 0.1); 
\draw[black](2.2,-1.2) circle(0.02);
\draw(2.1,-1.5) node{{\scriptsize $v_1$}};
\fill(2.35, -1.2) circle(0.5pt); 
\fill(2.5, -1.2) circle(0.5pt); 
\fill(2.65, -1.2) circle(0.5pt); 
\draw[black](2.95,-1.2) circle(0.02);
\draw(2.95,-1.5) node{{\scriptsize $v_2$}};
\fill(3.25, -1.2) circle(0.5pt); 
\fill(3.4, -1.2) circle(0.5pt); 
\fill(3.55, -1.2) circle(0.5pt); 
\fill(4.0, -1.2) circle(0.5pt); 
\fill(4.15, -1.2) circle(0.5pt); 
\fill(4.3, -1.2) circle(0.5pt); 
\draw[thick,black](4.9,-1.2) ellipse(0.5 and 0.1); 
\fill(4.75, -1.2) circle(0.5pt); 
\fill(4.9, -1.2) circle(0.5pt); 
\fill(5.05, -1.2) circle(0.5pt); 
\draw[black](5.3,-1.2) circle(0.02);
\draw(5.3,-1.5) node{{\scriptsize $v_{a+1}$}};
\draw(3.6,0) node{{\scriptsize $a$}};
\draw(3.6,-0.5) node{$\overbrace{~~~~~~~~~~~~~~~~~~~~~~~}$};
\draw[thick,black,rotate around={90:(5.3,-1.2)}](5.75,-1.2) ellipse(0.5 and 0.1);
\fill(5.3, -0.85) circle(0.5pt); 
\fill(5.3, -0.7) circle(0.5pt); 
\fill(5.3, -0.55) circle(0.5pt); 
\draw[black](5.3,-0.35) circle(0.02);
\draw[thick,black](5.7,-1.2) ellipse(0.5 and 0.1); 
\fill(5.6,-1.2) circle(0.5pt); 
\fill(5.75,-1.2) circle(0.5pt); 
\fill(5.9,-1.2) circle(0.5pt); 
\fill(6.35,-1.2) circle(0.5pt); 
\fill(6.5,-1.2) circle(0.5pt); 
\fill(6.65,-1.2) circle(0.5pt); 
\draw[thick,black](7.3,-1.2) ellipse(0.5 and 0.1); 
\fill(7.1,-1.2) circle(0.5pt); 
\fill(7.25,-1.2) circle(0.5pt); 
\fill(7.4,-1.2) circle(0.5pt); 
\draw[black](7.7,-1.2) circle(0.02);
\draw(7.7,-1.5) node{{\scriptsize $v_{a+b}$}};
\draw[thick,black](8.1,-1.2) ellipse(0.5 and 0.1); 
\fill(8.0,-1.2) circle(0.5pt); 
\fill(8.15,-1.2) circle(0.5pt); 
\fill(8.3,-1.2) circle(0.5pt); 
\draw[black](8.5,-1.2) circle(0.02);
\draw(8.9,-1.5) node{{\scriptsize $v_{a+b+1}$}};
\draw(7,0) node{{\scriptsize $b$}};
\draw(7,-0.5) node{$\overbrace{~~~~~~~~~~~~~~~~~~~~~~~}$};
\draw(5.3,-2) node{{\scriptsize $\mathcal{T}(a,b)$}};
\end{tikzpicture}
\begin{tikzpicture}
\draw[thick,black](2.6,-1.2) ellipse(0.5 and 0.1); 
\draw[thick,black](3.3,-1.2) ellipse(0.5 and 0.1); 
\draw[black](2.2,-1.2) circle(0.02);
\draw(2.1,-1.5) node{{\scriptsize $v_1$}};
\fill(2.35, -1.2) circle(0.5pt); 
\fill(2.5, -1.2) circle(0.5pt); 
\fill(2.65, -1.2) circle(0.5pt); 
\draw[black](2.95,-1.2) circle(0.02);
\draw(2.95,-1.5) node{{\scriptsize $v_2$}};
\fill(3.25, -1.2) circle(0.5pt); 
\fill(3.4, -1.2) circle(0.5pt); 
\fill(3.55, -1.2) circle(0.5pt); 
\fill(4.0, -1.2) circle(0.5pt); 
\fill(4.15, -1.2) circle(0.5pt); 
\fill(4.3, -1.2) circle(0.5pt); 
\draw[thick,black](4.9,-1.2) ellipse(0.5 and 0.1); 
\fill(4.75, -1.2) circle(0.5pt); 
\fill(4.9, -1.2) circle(0.5pt); 
\fill(5.05, -1.2) circle(0.5pt); 
\draw[black](5.3,-1.2) circle(0.02);
\draw(5.3,-1.5) node{{\scriptsize $v_{a+1}$}};
\draw(3.6,0) node{{\scriptsize $a$}};
\draw(3.6,-0.5) node{$\overbrace{~~~~~~~~~~~~~~~~~~~~~~~}$};
\draw[thick,black,rotate around={90:(5.3,-1.2)}](5.75,-1.2) ellipse(0.5 and 0.1);
\fill(5.3, -0.85) circle(0.5pt); 
\fill(5.3, -0.7) circle(0.5pt); 
\fill(5.3, -0.55) circle(0.5pt); 
\draw[black](5.3,-0.35) circle(0.02);
\draw[thick,black](5.7,-1.2) ellipse(0.5 and 0.1); 
\fill(5.6,-1.2) circle(0.5pt); 
\fill(5.75,-1.2) circle(0.5pt); 
\fill(5.9,-1.2) circle(0.5pt); 
\fill(6.35,-1.2) circle(0.5pt); 
\fill(6.5,-1.2) circle(0.5pt); 
\fill(6.65,-1.2) circle(0.5pt); 
\draw[thick,black](7.3,-1.2) ellipse(0.5 and 0.1); 
\fill(7.1,-1.2) circle(0.5pt); 
\fill(7.25,-1.2) circle(0.5pt); 
\fill(7.4,-1.2) circle(0.5pt); 
\draw[black](7.7,-1.2) circle(0.02);
\draw(7.7,-1.5) node{{\scriptsize $v_{a+b+1}$}};
\draw[thick,black,rotate around={90:(7.7,-1.2)}](8.15,-1.2) ellipse(0.5 and 0.1);
\fill(7.7, -0.85) circle(0.5pt); 
\fill(7.7, -0.7) circle(0.5pt); 
\fill(7.7, -0.55) circle(0.5pt); 
\draw[black](7.7,-0.35) circle(0.02);
\draw[thick,black](8.1,-1.2) ellipse(0.5 and 0.1); 
\fill(8.0,-1.2) circle(0.5pt); 
\fill(8.15,-1.2) circle(0.5pt); 
\fill(8.3,-1.2) circle(0.5pt); 
\draw(6.5,0) node{{\scriptsize $b$}};
\draw(6.5,-0.5) node{$\overbrace{~~~~~~~~~~~~~~}$};
\fill(8.75,-1.2) circle(0.5pt); 
\fill(8.9,-1.2) circle(0.5pt); 
\fill(9.05,-1.2) circle(0.5pt); 
\draw[thick,black](9.7,-1.2) ellipse(0.5 and 0.1); 
\fill(9.55,-1.2) circle(0.5pt); 
\fill(9.7,-1.2) circle(0.5pt); 
\fill(9.85,-1.2) circle(0.5pt); 
\draw[black](10.1,-1.2) circle(0.02);
\draw(10.1,-1.5) node{{\scriptsize $v_{a+b+c}$}};
\draw[thick,black](10.5,-1.2) ellipse(0.5 and 0.1); 
\fill(10.35,-1.2) circle(0.5pt); 
\fill(10.5,-1.2) circle(0.5pt); 
\fill(10.65,-1.2) circle(0.5pt); 
\draw[black](10.9,-1.2) circle(0.02);
\draw(11.5,-1.5) node{{\scriptsize $v_{a+b+c+1}$}};
\draw(9.35,0) node{{\scriptsize $c$}};
\draw(9.35,-0.5) node{$\overbrace{~~~~~~~~~~~~~~~~~~~~~~}$};
\draw(6.5,-2) node{{\scriptsize $\mathcal{Q}(a,b,c)$}};
\end{tikzpicture}
\begin{tikzpicture}
\draw[thick,black](2.6,-1.2) ellipse(0.5 and 0.1); 
\draw[thick,black](3.3,-1.2) ellipse(0.5 and 0.1); 
\draw[black](2.2,-1.2) circle(0.02);
\fill(2.35, -1.2) circle(0.5pt); 
\draw(2.1,-1.5) node{{\scriptsize $v_1$}};
\fill(2.5, -1.2) circle(0.5pt); 
\fill(2.65, -1.2) circle(0.5pt); 
\draw[black](2.95,-1.2) circle(0.02);
\draw(2.95,-1.5) node{{\scriptsize $v_2$}};
\fill(3.25, -1.2) circle(0.5pt); 
\fill(3.4, -1.2) circle(0.5pt); 
\fill(3.55, -1.2) circle(0.5pt); 
\fill(4.0, -1.2) circle(0.5pt); 
\fill(4.15, -1.2) circle(0.5pt); 
\fill(4.3, -1.2) circle(0.5pt); 
\draw[thick,black](4.9,-1.2) ellipse(0.5 and 0.1); 
\fill(4.75, -1.2) circle(0.5pt); 
\fill(4.9, -1.2) circle(0.5pt); 
\fill(5.05, -1.2) circle(0.5pt); 
\draw[black](5.3,-1.2) circle(0.02);
\draw(5.3,-1.5) node{{\scriptsize $v_{a+1}$}};
\draw(3.6,0) node{{\scriptsize $a$}};
\draw(3.6,-0.5) node{$\overbrace{~~~~~~~~~~~~~~~~~~~~~~~}$};
\draw[thick,black,rotate around={90:(5.3,-1.2)}](5.75,-1.2) ellipse(0.5 and 0.1);
\fill(5.3, -0.85) circle(0.5pt); 
\fill(5.3, -0.7) circle(0.5pt); 
\fill(5.3, -0.55) circle(0.5pt); 
\draw[black](5.3,-0.35) circle(0.02);
\draw[thick,black](5.70,-1.2) ellipse(0.5 and 0.1); 
\fill(5.6,-1.2) circle(0.5pt); 
\fill(5.75,-1.2) circle(0.5pt); 
\fill(5.9,-1.2) circle(0.5pt); 
\fill(6.35,-1.2) circle(0.5pt); 
\fill(6.5,-1.2) circle(0.5pt); 
\fill(6.65,-1.2) circle(0.5pt); 
\draw[thick,black](7.3,-1.2) ellipse(0.5 and 0.1); 
\fill(7.1,-1.2) circle(0.5pt); 
\fill(7.25,-1.2) circle(0.5pt); 
\fill(7.4,-1.2) circle(0.5pt); 
\draw[black](7.7,-1.2) circle(0.02);
\draw(7.7,-1.5) node{{\scriptsize $v_{a+b+1}$}};
\draw[thick,black,rotate around={90:(7.7,-1.2)}](8.15,-1.2) ellipse(0.5 and 0.1);
\fill(7.7, -0.85) circle(0.5pt); 
\fill(7.7, -0.7) circle(0.5pt); 
\fill(7.7, -0.55) circle(0.5pt); 
\draw[black](7.7,-0.35) circle(0.02);
\draw[thick,black](8.10,-1.2) ellipse(0.5 and 0.1); 
\fill(8.0,-1.2) circle(0.5pt); 
\fill(8.15,-1.2) circle(0.5pt); 
\fill(8.3,-1.2) circle(0.5pt); 
\draw(6.5,0) node{{\scriptsize $b$}};
\draw(6.5,-0.5) node{$\overbrace{~~~~~~~~~~~~~~}$};
\fill(8.75,-1.2) circle(0.5pt); 
\fill(8.9,-1.2) circle(0.5pt); 
\fill(9.05,-1.2) circle(0.5pt); 
\draw[thick,black](9.7,-1.2) ellipse(0.5 and 0.1); 
\fill(9.55,-1.2) circle(0.5pt); 
\fill(9.7,-1.2) circle(0.5pt); 
\fill(9.85,-1.2) circle(0.5pt); 
\draw[black](10.1,-1.2) circle(0.02);
\draw(10.1,-1.5) node{{\scriptsize $v_{a+b+c+1}$}};
\draw[thick,black,rotate around={90:(10.1,-1.2)}](10.55,-1.2) ellipse(0.5 and 0.1);
\fill(10.1, -0.85) circle(0.5pt); 
\fill(10.1, -0.7) circle(0.5pt); 
\fill(10.1, -0.55) circle(0.5pt); 
\draw[black](10.1,-0.35) circle(0.02);
\draw(8.9,0) node{{\scriptsize $c$}};
\draw(8.9,-0.5) node{$\overbrace{~~~~~~~~~~~~~~~~}$};
\draw[thick,black](10.5,-1.2) ellipse(0.5 and 0.1); 
\fill(10.35,-1.2) circle(0.5pt); 
\fill(10.5,-1.2) circle(0.5pt); 
\fill(10.65,-1.2) circle(0.5pt); 
\fill(11.1,-1.2) circle(0.5pt); 
\fill(11.25,-1.2) circle(0.5pt); 
\fill(11.4,-1.2) circle(0.5pt); 
\draw[thick,black](12.0,-1.2) ellipse(0.5 and 0.1); 
\fill(11.85,-1.2) circle(0.5pt); 
\fill(12.0,-1.2) circle(0.5pt); 
\fill(12.15,-1.2) circle(0.5pt); 
\draw[black](12.4,-1.2) circle(0.02);
 \draw(12.4,-1.5) node{{\scriptsize $v_{a+b+c+d}$}};
\draw[thick,black](12.8,-1.2) ellipse(0.5 and 0.1); 
\fill(12.65,-1.2) circle(0.5pt); 
\fill(12.8,-1.2) circle(0.5pt); 
\fill(12.95,-1.2) circle(0.5pt); 
\draw[black](13.2,-1.2) circle(0.02);
\draw(14,-1.5) node{{\scriptsize $v_{a+b+c+d+1}$}};
\draw(11.7,0) node{{\scriptsize $d$}};
\draw(11.7,-0.5) node{$\overbrace{~~~~~~~~~~~~~~~~~~~~~}$};
\draw(7.7,-2) node{{\scriptsize $\mathcal{R}(a,b,c,d)$}};
\end{tikzpicture}
\caption{\label{figTQR} $\mathcal{T}(a,b)$,~$\mathcal{Q}(a,b,c)$~and~$\mathcal{R}(a,b,c,d)$}
\end{figure}
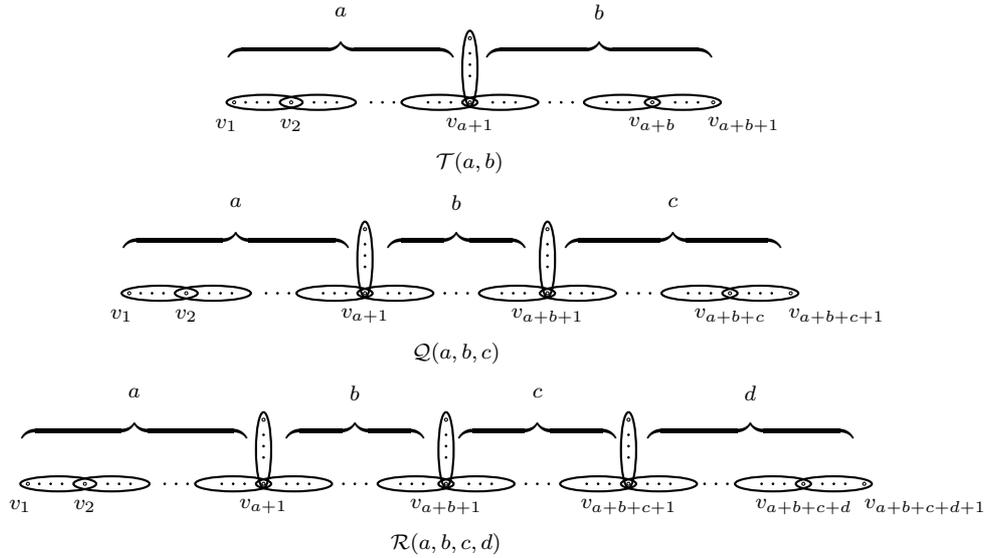

For two given $r$-uniform supertrees, Theorem \ref{wuxian} can provide us with a simple method to  investigate whether their SR \textcolor{black}{and ME are the same or not}.
   Next, we give an example to show how to apply \textcolor{black}{Lemmas  \ref{large} and \ref{th1}} to determine
   that the SR \textcolor{black}{and ME} of two $r$-uniform  supertrees are the same, which is shown in  Theorem \ref{wth4}.

  Let $H$ be an ordinary graph.
  The  $r$-th power of $H$,  denoted by $\mathcal{H}^r$, is obtained from $H$ by adding $(r-2)$ new vertices into  each edge of $H$, where $r\geq 3$.
  Let $P_t$ be a path of length $t$ and $\mathcal{P}^r_t$ be its $r$-th power, where $t\geq 0$.
  When $t=0$, $P_0$ is a vertex.
  We call $\mathcal{P}^r_t$ a loose path of length $t$.
   Let $\mathcal{P}_{t}^r=v_1e_1v_2e_2v_3\cdots v_{t}e_{t}v_{t+1}$, where $r\geq 2$, $t\geq 1$ and $e_{i}=\{v_i, u_{i,1}, \cdots, u_{i,r-2},v_{i+1}\}$ with $i=1,2,\cdots,t$.
  Let  $\mathcal{T}(a,b)$, $\mathcal{Q}(a,b,c)$ and  $\mathcal{R}(a,b,c,d)$ be  three hypergraphs defined as follows, where $a$, $b$ and $c$ are three positive integers.     $\mathcal{T}(a,b)$ is obtained from  $\mathcal{P}^r_{a+b}$ by  attaching one pendent edge at vertex $v_{a+1}$ of $\mathcal{P}^r_{a+b}$,
  $\mathcal{Q}(a,b,c)$ is  obtained from  $\mathcal{P}^r_{a+b+c}$ by attaching one pendent edge at vertices $v_{a+1}$ and $v_{a+b+1}$  of  $\mathcal{P}^r_{a+b+c}$,
 and $\mathcal{R}(a,b,c,d)$ is obtained from  $\mathcal{P}^r_{a+b+c+d}$ by attaching one pendent edge at vertices  $v_{a+1}$,  $v_{a+b+1}$ and $v_{a+b+c+1}$ of
 $\mathcal{P}^r_{a+b+c+d}$.
 $\mathcal{T}(a,b)$, $\mathcal{Q}(a,b,c)$ and  $\mathcal{R}(a,b,c,d)$ are shown in Fig.\ \ref{figTQR}.
 In particular, in  $\mathcal{P}_{t}^r$, if $r=2$, then  $\mathcal{P}_{t}^r$ is the  path $P_t$.
 Furthermore, when $r=2$,  $\mathcal{T}(a,b)$, $\mathcal{Q}(a,b,c)$ and  $\mathcal{R}(a,b,c,d)$ are graphs and written as  ${T}(a,b)$, $Q(a,b,c)$ and  $R(a,b,c,d)$, respectively.

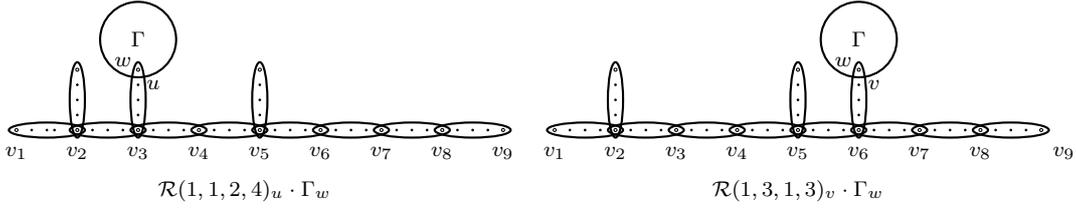
\begin{figure}
\thicklines
\begin{tikzpicture}
\thicklines
\draw[thick,black](2.6,-1.2) ellipse(0.5 and 0.1); 
\draw[black](2.2,-1.2) circle(0.02);
\draw(2.2,-1.5) node{{\scriptsize $v_1$}};
\fill(2.4, -1.2) circle(0.5pt); 
\fill(2.6, -1.2) circle(0.5pt); 
\fill(2.7, -1.2) circle(0.5pt); 
\draw[black](3.0,-1.2) circle(0.02);
\draw(3.0,-1.5) node{{\scriptsize $v_2$}};
\draw[thick,black,rotate around={90:(3.0,-1.2)}](3.4,-1.2) ellipse(0.5 and 0.1);
\fill(3.0, -0.6) circle(0.5pt); 
\fill(3.0, -0.8) circle(0.5pt); 
\fill(3.0, -1.0) circle(0.5pt); 
\draw[black](3.0,-0.4) circle(0.02);
\draw[thick,black](3.4,-1.2) ellipse(0.5 and 0.1); 
\fill(3.2, -1.2) circle(0.5pt); 
\fill(3.4, -1.2) circle(0.5pt); 
\fill(3.6, -1.2) circle(0.5pt); 
\draw[black](3.8,-1.2) circle(0.02);
\draw(3.8,-1.5) node{{\scriptsize $v_3$}};
\draw[thick,black,rotate around={90:(3.8,-1.2)}](4.2,-1.2) ellipse(0.5 and 0.1);
\fill(3.8, -0.6) circle(0.5pt); 
\fill(3.8, -0.8) circle(0.5pt); 
\fill(3.8, -1.0) circle(0.5pt); 
\draw[black](3.8,-0.4) circle(0.02);
\draw(4,-0.6) node{{\scriptsize $u$}};
\draw(3.6,-0.3) node{{\scriptsize $w$}};
\draw[thick,black](3.8,0) circle(0.5);
\draw(3.8,0) node{{\scriptsize $\Gamma$}};
\draw[thick,black](4.2,-1.2) ellipse(0.5 and 0.1); 
\fill(4.0, -1.2) circle(0.5pt); 
\fill(4.2, -1.2) circle(0.5pt); 
\fill(4.4, -1.2) circle(0.5pt); 
\draw[black](4.6,-1.2) circle(0.02);
\draw(4.6,-1.5) node{{\scriptsize $v_4$}};
\draw[thick,black](5,-1.2) ellipse(0.5 and 0.1); 
\fill(4.8, -1.2) circle(0.5pt); 
\fill(5, -1.2) circle(0.5pt); 
\fill(5.2, -1.2) circle(0.5pt); 
\draw[black](5.4,-1.2) circle(0.02);
\draw(5.4,-1.5) node{{\scriptsize $v_5$}};
\draw[thick,black,rotate around={90:(5.4,-1.2)}](5.8,-1.2) ellipse(0.5 and 0.1);
\fill(5.4, -0.6) circle(0.5pt); 
\fill(5.4, -0.8) circle(0.5pt); 
\fill(5.4, -1.0) circle(0.5pt); 
\draw[black](5.4,-0.4) circle(0.02);
\draw[thick,black](5.8,-1.2) ellipse(0.5 and 0.1); 
\fill(5.6, -1.2) circle(0.5pt); 
\fill(5.8, -1.2) circle(0.5pt); 
\fill(6.0, -1.2) circle(0.5pt); 
\draw[black](6.2,-1.2) circle(0.02);
\draw(6.2,-1.5) node{{\scriptsize $v_6$}};
\draw[thick,black](6.6,-1.2) ellipse(0.5 and 0.1); 
\fill(6.4, -1.2) circle(0.5pt); 
\fill(6.6, -1.2) circle(0.5pt); 
\fill(6.8, -1.2) circle(0.5pt); 
\draw[black](7.0,-1.2) circle(0.02);
\draw(7.0,-1.5) node{{\scriptsize $v_7$}};
\draw[thick,black](7.4,-1.2) ellipse(0.5 and 0.1); 
\fill(7.2, -1.2) circle(0.5pt); 
\fill(7.4, -1.2) circle(0.5pt); 
\fill(7.6, -1.2) circle(0.5pt); 
\draw[black](7.8,-1.2) circle(0.02);
\draw(7.8,-1.5) node{{\scriptsize $v_8$}};
\draw[thick,black](8.2,-1.2) ellipse(0.5 and 0.1); 
\fill(8.0, -1.2) circle(0.5pt); 
\fill(8.2, -1.2) circle(0.5pt); 
\fill(8.4, -1.2) circle(0.5pt); 
\draw[black](8.6,-1.2) circle(0.02);
\draw(8.6,-1.5) node{{\scriptsize $v_9$}};
\draw(5.2,-2) node{{\scriptsize $\mathcal{R}(1,1,2,4)_u\cdot\Gamma_w$}};
\end{tikzpicture}
\begin{tikzpicture}
\draw[thick,black](2.6,-1.2) ellipse(0.5 and 0.1); 
\draw[black](2.2,-1.2) circle(0.02);
\draw(2.2,-1.5) node{{\scriptsize $v_1$}};
\fill(2.4, -1.2) circle(0.5pt); 
\fill(2.6, -1.2) circle(0.5pt); 
\fill(2.8, -1.2) circle(0.5pt); 
\draw[black](3.0,-1.2) circle(0.02);
\draw(3.0,-1.5) node{{\scriptsize $v_2$}};
\draw[thick,black,rotate around={90:(3.0,-1.2)}](3.4,-1.2) ellipse(0.5 and 0.1);
\fill(3.0, -0.6) circle(0.5pt); 
\fill(3.0, -0.8) circle(0.5pt); 
\fill(3.0, -1) circle(0.5pt); 
\draw[black](3.0,-0.4) circle(0.02);
\draw[thick,black](3.4,-1.2) ellipse(0.5 and 0.1); 
\fill(3.2, -1.2) circle(0.5pt); 
\fill(3.4, -1.2) circle(0.5pt); 
\fill(3.6, -1.2) circle(0.5pt); 
\draw[black](3.8,-1.2) circle(0.02);
\draw(3.8,-1.5) node{{\scriptsize $v_3$}};
\draw[thick,black](4.2,-1.2) ellipse(0.5 and 0.1); 
\fill(4, -1.2) circle(0.5pt); 
\fill(4.2, -1.2) circle(0.5pt); 
\fill(4.4, -1.2) circle(0.5pt); 
\draw[black](4.6,-1.2) circle(0.02);
\draw(4.6,-1.5) node{{\scriptsize $v_4$}};
\draw[thick,black](5,-1.2) ellipse(0.5 and 0.1); 
\fill(4.8, -1.2) circle(0.5pt); 
\fill(5, -1.2) circle(0.5pt); 
\fill(5.2, -1.2) circle(0.5pt); 
\draw[black](5.4,-1.2) circle(0.02);
\draw(5.4,-1.5) node{{\scriptsize $v_5$}};
\draw[thick,black,rotate around={90:(5.4,-1.2)}](5.8,-1.2) ellipse(0.5 and 0.1);
\fill(5.4, -0.6) circle(0.5pt); 
\fill(5.4, -0.8) circle(0.5pt); 
\fill(5.4, -1.0) circle(0.5pt); 
\draw[black](5.4,-0.4) circle(0.02);
\draw[thick,black](5.8,-1.2) ellipse(0.5 and 0.1); 
\fill(5.6, -1.2) circle(0.5pt); 
\fill(5.8, -1.2) circle(0.5pt); 
\fill(6.0, -1.2) circle(0.5pt); 
\draw[black](6.2,-1.2) circle(0.02);
\draw(6.2,-1.5) node{{\scriptsize $v_6$}};
\draw[thick,black,rotate around={90:(6.2,-1.2)}](6.6,-1.2) ellipse(0.5 and 0.1);
\fill(6.2, -0.6) circle(0.5pt); 
\fill(6.2, -0.8) circle(0.5pt); 
\fill(6.2, -1.0) circle(0.5pt); 
\draw[black](6.2,-0.4) circle(0.02);
\draw(6.4,-0.6) node{{\scriptsize $v$}};
\draw(6.0,-0.3) node{{\scriptsize $w$}};
\draw[thick,black](6.2,0) circle(0.5);
\draw(6.2,0) node{{\scriptsize $\Gamma$}};
\draw[thick,black](6.6,-1.2) ellipse(0.5 and 0.1); 
\fill(6.4, -1.2) circle(0.5pt); 
\fill(6.6, -1.2) circle(0.5pt); 
\fill(6.8, -1.2) circle(0.5pt); 
\draw[black](7.0,-1.2) circle(0.02);
\draw(7.0,-1.5) node{{\scriptsize $v_7$}};
\draw[thick,black](7.4,-1.2) ellipse(0.5 and 0.1); 
\fill(7.2, -1.2) circle(0.5pt); 
\fill(7.4, -1.2) circle(0.5pt); 
\fill(7.6, -1.2) circle(0.5pt); 
\draw[black](7.8,-1.2) circle(0.02);
\draw(7.8,-1.5) node{{\scriptsize $v_8$}};
\draw[thick,black](8.2,-1.2) ellipse(0.5 and 0.1); 
\fill(8.0, -1.2) circle(0.5pt); 
\fill(8.2, -1.2) circle(0.5pt); 
\fill(8.4, -1.2) circle(0.5pt); 
\draw[black](8.6,-1.2) circle(0.02);
\draw(8.9,-1.5) node{{\scriptsize $v_9$}};
\draw(5.4,-2) node{{\scriptsize $\mathcal{R}(1,3,1,3)_v\cdot\Gamma_w$}};
\end{tikzpicture}
\caption{\label{eg} $\mathcal{R}(1,1,2,4)_u\cdot\Gamma_w$~and~$\mathcal{R}(1,3,1,3)_v\cdot\Gamma_w$}
\end{figure}

 \begin{theorem} \label{wth4}
  Let  $\Gamma$ be  an $r$-uniform supertree with $w\in V(\Gamma)$, where $r\geq 3$.
  We have $\rho(\mathcal{R}(1,1,2,4)_u \cdot \Gamma_w)=\rho(\mathcal{R}(1,3,1,3)_v \cdot\Gamma_w)$ and \textcolor{black}{$ME(\mathcal{R}(1,1,2,4)_u \cdot \Gamma_w)=ME(\mathcal{R}(1,3,1,3)_v \cdot\Gamma_w)$}, where
  $u$ of $\mathcal{R}(1,1,2,4)$ (respectively, $v$ of $\mathcal{R}(1,3,1,3)$)
  is a pendent vertex of the pendent edge attached at $v_{3}$ of $\mathcal{P}^r_{8}$ of $\mathcal{R}(1,1,2,4)$
  (respectively, $v_{6}$ of $\mathcal{P}^r_{8}$ of $\mathcal{R}(1,3,1,3)$).
   $\mathcal{R}(1,1,2,4)_u\cdot \Gamma_w$ and $\mathcal{R}(1,3,1,3)_v\cdot \Gamma_w$ are shown in Fig.\ \ref{eg}.
  \end{theorem}

\noindent\textbf{Proof}. For simplicity, let $\mathcal{R}(1,1,2,4)=\mathcal{G'}$ and $\mathcal{R}(1,3,1,3)=\mathcal{H'}$.
   We have
  \begin{align}
  \varphi(\mathcal{G'},x)&=x\varphi(\mathcal{G'}-u,x)-\sum_{e\in E_{\mathcal{G}}(u)}
  \varphi(\mathcal{G'}-V(e),x)\label{1.1}\\
  &=x\varphi\left[\mathcal{Q}(1,3,4)\cup N_{r-2},x\right]-\varphi\left[\mathcal{P}_{2}^r\cup \mathcal{T}(1,4)\cup N_{2(r-2)},x\right]\label{1.2}\\
  &=x^{r-1}\varphi(\mathcal{Q}(1,3,4),x)-x^{2(r-2)}\varphi(\mathcal{P}_{2}^r,x)\varphi(\mathcal{T}(1,4),x).\label{1.3}
 \end{align}
 It is noted that  \eqref{1.1}  follows from Lemma \ref{su}(b),
 \eqref{1.2} is derived from $\mathcal{G'}-u\cong \mathcal{Q}(1,3,4)\cup N_{r-2}$ and
 $\mathcal{G'}-V(e)\cong \mathcal{P}_{2}^r \cup \mathcal{T}(1,4)\cup N_{2(r-2)}$, and
  \eqref{1.3}  is deduced from Lemma \ref{su}(a).
 Similarly, we obtain  that the expression of $\varphi(\mathcal{H'},x)$ is the same as the right-hand  side of (\ref{1.3}).
  Namely, we obtain $\varphi(\mathcal{G'},x)=\varphi(\mathcal{H'},x)$.
  Obviously, we can check that $\mathcal{G'}-u \cong \mathcal{H'}-v$.
    \textcolor{black}{Thus, by Lemma \ref{th1}, we get $\varphi(\mathcal{R}(1,1,2,4)_u \cdot \Gamma_w)=\varphi(\mathcal{R}(1,3,1,3)_v \cdot\Gamma_w)$.
    Furthermore, by Lemma \ref{large},  we obtain $\rho(\mathcal{R}(1,1,2,4)_u \cdot \Gamma_w)=\rho(\mathcal{R}(1,3,1,3)_v \cdot\Gamma_w)$}.
    \textcolor{black}{By the definition of the matching energy of an $r$-uniform  hypergraph, we obtain $ME(\mathcal{R}(1,1,2,4)_u \cdot \Gamma_w)=ME(\mathcal{R}(1,3,1,3)_v \cdot\Gamma_w)$.}
      ~~$\Box$

  From Theorem \ref{wth4}, we know that the problem of  determining whether  two $r$-uniform supertrees have the same SR \textcolor{black}{and ME} can be converted into  the problem of  investigating the properties of their subgraphs.

  For graphs, by  Lemmas \ref{tp} and \ref{th1}, we have Theorem \ref{c1} as follows.
  It should be noted that \textcolor{black}{  Theorem \ref{c1}(i)} is  a natural generalization of  Lemma \ref{th1} when $r=2$ and it can be found on Page 159  in Ref. \cite{Cvetkovi1980}.

  \begin{theorem}\label{c1}
   Let $G$, $H$ and $\Gamma$ be three trees, where $G$ and $H$ have the same number of vertices.
   Let $u\in V(G)$ and $v\in V(H)$.
   If $\phi(G,x)=\phi(H,x)$
   and $\phi(G-u,x)=\phi(H-v,x)$, then for any $w\in V(\Gamma)$,
   we have \textcolor{black}{ (i). $\phi(G(u,w)\Gamma,x)=\phi(H(v,w)\Gamma,x)$.
   Namely,  $G(u,w)\Gamma$ and $H(v,w)\Gamma$ are $\bm{M}$-cospectral, where  $\bm{M}$ is the adjacency matrix.
   (ii). $ME(G(u,w)\Gamma,x)=ME(H(v,w)\Gamma,x)$.}
   \end{theorem}

  By Theorem \ref{c1} and the  methods similar to those for Theorem \ref{c2w}, we get Theorem \ref{c11} as follows.
  \textcolor{black}{It is noted that Theorem \ref{c11}(i)} can be found
 on Page 158  in Ref. \cite{Schwenk1974}.

 \begin{theorem} \label{c11}
 Let $G$ and $H$ be two trees with the same number of vertices.
 Let  $u\in V(G)$,  $v\in V(H)$ and $m$ be a positive integer.
 If $\phi(G,x)=\phi(H,x)$ and $\phi(G-u,x)=\phi(H-v,x)$,
 then we have
 (i). $\phi(H_v^m)=\phi(H_v^{m-1}\cdot G_u)
 =\phi(H^{m-2}_v\cdot G_u^2)=\cdots
 =\phi(H_v\cdot G_u^{m-1})
 =\phi(G_u^m)$;
 (ii). \textcolor{black}{$ME(H_v^m)=ME(H_v^{m-1}\cdot G_u)
 =ME(H^{m-2}_v\cdot G_u^2)=\cdots
 =ME(H_v\cdot G_u^{m-1})
 =ME(G_u^m)$.}
\end{theorem}

\section{The second pair of  $r$-uniform supertrees with the same spectral radius and matching energy}

  In this section, we construct the second   pair of $r$-uniform supertrees \textcolor{black}{with the same SR and ME}, which is shown in  Theorem \ref{c4}, where $r\geq 3$.
  In Theorem \ref{c4}. Since $m$ is a variable,  an infinite families of  $r$-uniform supertrees with the same SR \textcolor{black}{and  ME} are also deduced.
 To obtain our results, Lemma  \ref{th2}  is introduced first.
 It is pointed out that Lemma  \ref{th2} generalizes many  known results in  the previous literatures.
 A  pair of graphs which is  $\bm{M}$-cospectral is deduced from Lemma  \ref{th2} (as shown in Theorem \ref{c3s}), where  $\bm{M}$ is the adjacency matrix.

%

\begin{figure}
\thicklines
\begin{tikzpicture}
\draw[thick,black] (0,-0.05) circle (0.5);
\draw(0,-0.05) node{{\scriptsize $\mathcal{G}$}};
\draw[thick,black,rotate around={-90:(0,-0.4)}](1,-0.4) ellipse(1.05 and 0.15);
\draw[black](0.015,-0.47) circle(0.02);
\draw(0.2,-0.3) node{{\scriptsize $u$}};
\fill(0, -0.9) circle(0.5pt); 
\fill(0, -1.2) circle(0.5pt); 
\fill(0, -1.5) circle(0.5pt); 
\draw[black](0,-2.35) circle(0.02);
\draw(0.2,-2.55) node{{\scriptsize $v$}};
\draw[thick,black](0,-2.8) circle(0.5);
\draw(0,-2.8) node{{\scriptsize $\mathcal{H}$}};
\draw[thick,black](-1.5,-2.8) circle(0.5);
\draw(-1.5,-2.8) node{{\scriptsize $\mathcal{H}$}};
\draw[thick,black,rotate around={-125:(0,-0.4)}](1.3,-0.35) ellipse(1.35 and 0.15);
\fill(-0.35, -0.9) circle(0.5pt); 
\fill(-0.55, -1.2) circle(0.5pt); 
\fill(-0.75, -1.5) circle(0.5pt); 
\draw[black](-1.35,-2.4) circle(0.02);
\draw(-1.65,-2.5) node{{\scriptsize $v$}};
\fill(2.5, -2.8) circle(0.5pt); 
\fill(2.0, -2.8) circle(0.5pt); 
\fill(1.5, -2.8) circle(0.5pt); 
\draw[thick,black](3.5,-2.8) circle(0.5);
\draw(3.5,-2.8) node{{\scriptsize $\mathcal{H}$}};
\draw[thick,black,rotate around={-30:(0,-0.4)}](1.95,-0.455) ellipse(2.0 and 0.15);
\fill(0.8, -0.9) circle(0.5pt); 
\fill(1.3, -1.2) circle(0.5pt); 
\fill(1.8, -1.5) circle(0.5pt); 
\draw[black](3.3,-2.37) circle(0.02);
\draw(3.6,-2.5) node{{\scriptsize $v$}};
\draw(1.7,-3.5) node{$\underbrace{~~~~~~~~~~~~~~~~~~~~~~~~}$};
\draw(1.7,-4.0) node{{\scriptsize $m-1$}};
\draw[thick,black](1.5,-0.05) circle(0.5);
\draw(1.5,-0.05) node{{\scriptsize $\mathcal{G}$}};
\fill(2.5, 0) circle(0.5pt); 
\fill(3, 0) circle(0.5pt); 
\fill(3.5, 0) circle(0.5pt); 
\draw[thick,black](4.5,-0.05) circle(0.5);
\draw(4.5,-0.05) node{{\scriptsize $\mathcal{G}$}};
\draw(3,1.3) node{{\scriptsize $m-1$}};
\draw(3,0.8) node{$\overbrace{~~~~~~~~~~~~~~~~~~~~~}$};
\draw(1.7,-4.5) node{{\scriptsize $(\mathcal{G}_u\cdot m\mathcal{H}_v)\cup (m-1)\mathcal{G}$}};
\end{tikzpicture}
\begin{tikzpicture}
\draw[thick,black] (0,-0.05) circle (0.5);
\draw(0,-0.05) node{{\scriptsize $\mathcal{H}$}};
\draw[thick,black,rotate around={-90:(0,-0.4)}](1,-0.4) ellipse(1.05 and 0.15);
\draw[black](0.015,-0.47) circle(0.02);
\draw(0.2,-0.3) node{{\scriptsize $v$}};
\fill(0, -0.9) circle(0.5pt); 
\fill(0, -1.2) circle(0.5pt); 
\fill(0, -1.5) circle(0.5pt); 
\draw[black](0,-2.35) circle(0.02);
\draw(0.2,-2.55) node{{\scriptsize $u$}};
\draw[thick,black](0,-2.8) circle(0.5);
\draw(0,-2.8) node{{\scriptsize $\mathcal{G}$}};
\draw[thick,black](-1.5,-2.8) circle(0.5);
\draw(-1.5,-2.8) node{{\scriptsize $\mathcal{G}$}};
\draw[thick,black,rotate around={-125:(0,-0.4)}](1.3,-0.35) ellipse(1.35 and 0.15);
\fill(-0.35, -0.9) circle(0.5pt); 
\fill(-0.55, -1.2) circle(0.5pt); 
\fill(-0.75, -1.5) circle(0.5pt); 
\draw[black](-1.35,-2.4) circle(0.02);
\draw(-1.65,-2.5) node{{\scriptsize $u$}};
\fill(2.5, -2.8) circle(0.5pt); 
\fill(2.0, -2.8) circle(0.5pt); 
\fill(1.5, -2.8) circle(0.5pt); 
\draw[thick,black](3.5,-2.8) circle(0.5);
\draw(3.5,-2.8) node{{\scriptsize $\mathcal{G}$}};
\draw[thick,black,rotate around={-30:(0,-0.4)}](1.95,-0.455) ellipse(2.0 and 0.15);
\fill(0.8, -0.9) circle(0.5pt); 
\fill(1.3, -1.2) circle(0.5pt); 
\fill(1.8, -1.5) circle(0.5pt); 
\draw[black](3.3,-2.37) circle(0.02);
\draw(3.6,-2.5) node{{\scriptsize $u$}};
\draw(1.7,-3.5) node{$\underbrace{~~~~~~~~~~~~~~~~~~~~~~~~}$};
\draw(1.7,-4.0) node{{\scriptsize $m-1$}};
\draw[thick,black](1.5,-0.05) circle(0.5);
\draw(1.5,-0.05) node{{\scriptsize $\mathcal{H}$}};
\fill(2.5, 0) circle(0.5pt); 
\fill(3, 0) circle(0.5pt); 
\fill(3.5, 0) circle(0.5pt); 
\draw[thick,black](4.5,-0.05) circle(0.5);
\draw(4.5,-0.05) node{{\scriptsize $\mathcal{H}$}};
\draw(3,1.3) node{{\scriptsize $m-1$}};
\draw(3,0.8) node{$\overbrace{~~~~~~~~~~~~~~~~~~~~~}$};
\draw(1.7,-4.5) node{{\scriptsize $(\mathcal{H}_v\cdot m\mathcal{G}_u)\cup (m-1)\mathcal{H}$}};
\end{tikzpicture}
\caption{\label{c2} $(\mathcal{G}_u\cdot m\mathcal{H}_v)\cup (m-1)\mathcal{G}$~and~$(\mathcal{H}_v\cdot m\mathcal{G}_u)\cup (m-1)\mathcal{H}$}
\end{figure}
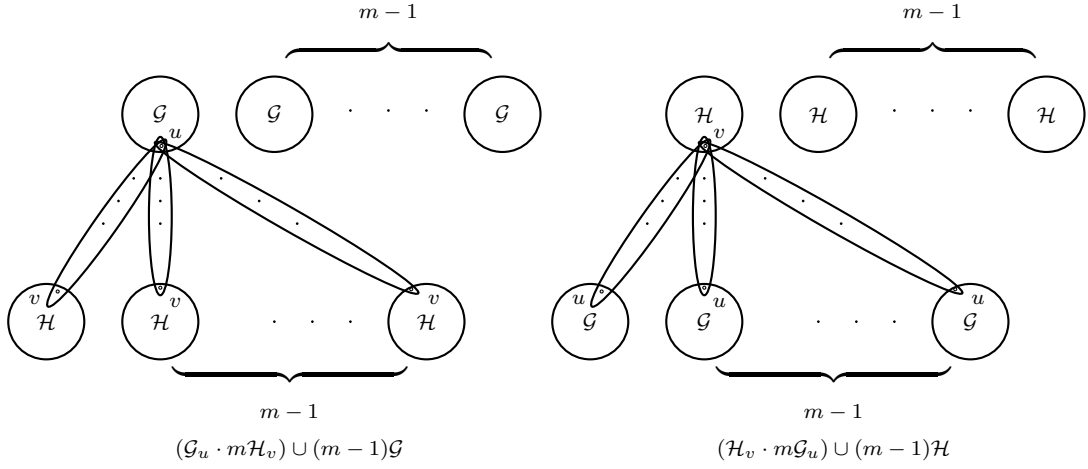

 Let $\mathcal{G}$ and  $\mathcal{H}$ be two $r$-uniform supertrees with $u\in V(\mathcal{G})$ and $v \in V(\mathcal{H})$, where $r\geq 2$.
   Let  $m$ be a  positive integer and $e_i=\{v_{i,1},\cdots, v_{i,r}\}$, where $i=1,\cdots,m$.
  Let $\mathcal{G}_u \cdot m\mathcal{H}_v$
  be the hypergraph obtained from $\mathcal{G}$, $m\mathcal{H}$ and $e_1, \cdots, e_m$  by identifying $v_{i,1}$ ($i=1,\cdots,m$) with
     $u$ of $\mathcal{G}$ and identifying $v_{i,r}$ ($i=1,\cdots,m$) with $v$ of each $\mathcal{H}$ of $m\mathcal{H}$ such that
   $\mathcal{G}_u \cdot m\mathcal{H}_v$ is also an $r$-uniform supertree.
   $\mathcal{G}_u \cdot m\mathcal{H}_v \cup (m-1)\mathcal{G}$ and
   $(\mathcal{H}_v\cdot m\mathcal{G}_u)\cup (m-1)\mathcal{H}$  are shown in Fig.\ \ref{c2}.
  Obviously, when $m=1$,  $\mathcal{G}_u \cdot m\mathcal{H}_v \cup (m-1)\mathcal{G} \cong (\mathcal{H}_v\cdot m\mathcal{G}_u)\cup (m-1)\mathcal{H}$.
  When
  $r=2$,  $\mathcal{G}_u \cdot m\mathcal{H}_v \cup (m-1)\mathcal{G}$ and
   $(\mathcal{H}_v\cdot m\mathcal{G}_u)\cup (m-1)\mathcal{H}$  are graphs and are written as
    $G_u \cdot m H_v \cup (m-1)G$ and
   $(H_v\cdot m G_u)\cup (m-1)H$, respectively.


 \begin{lemma}\label{th2}

  Let $\mathcal{G}$ and $\mathcal{H}$  be two $r$-uniform hypergraphs with $u\in V(\mathcal{G})$ and $v\in V(\mathcal{H})$, where  $r\geq 2$.
  We have     $\varphi\left[\left(\mathcal{G}_u\cdot m\mathcal{H}_v\right)\cup (m-1)\mathcal{G},x\right]=\varphi\left[\left(\mathcal{H}_v\cdot m\mathcal{G}_u\right)\cup (m-1)\mathcal{H},x\right]$, where $m$ is a positive integer.
 \end{lemma}

\noindent\textbf{Proof}. 
 Since $\mathcal{G}_u\cdot m\mathcal{H}_v-u\cong (\mathcal{G}-u)\cup m\mathcal{H}\cup N_{m(r-2)}$ and
 $E_{\mathcal{G}_u\cdot m\mathcal{H}_v}(u)=E_{\mathcal{G}}(u)\cup \{e_1,\cdots, e_m\}$,
  by Lemma \ref{su}(b), we obtain
  \begin{align}
    &\varphi(\mathcal{G}_u\cdot m\mathcal{H}_v,x)=x\varphi\left[\left(\mathcal{G}_u\cdot m\mathcal{H}_v\right)-u,x\right]-\sum_{e\in E_{\mathcal{G}_u\cdot m\mathcal{H}_v}(u)}\varphi\left[\left(\mathcal{G}_u\cdot m\mathcal{H}_v\right)-V(e),x\right]\label{2.1}\\
    &=x\varphi\left[\left(\mathcal{G}-u\right)\cup m\mathcal{H}\cup N_{m(r-2)},x\right]-\sum_{e\in E_{\mathcal{G}}(u)}\varphi\left[\left(\mathcal{G}-V(e)\right)\cup m\mathcal{H}\cup N_{m(r-2)},x\right]\nonumber\\
    &\quad-m\varphi\left[\left(\mathcal{G}-u\right)\cup \left(\mathcal{H}-v\right)\cup \left(m-1\right)\mathcal{H}\cup N_{(m-1)(r-2)},x\right].\label{2.2}
      \end{align}
  By   Lemma \ref{su}(a) and  extracting the common factors $x^{m(r-2)}$ and $\varphi(m\mathcal{H},x)=\varphi^m(\mathcal{H},x)$  from the first and the second terms on the right-hand side  of (\ref{2.2}), we obtain
    \begin{align}
    \varphi(\mathcal{G}_u\cdot m\mathcal{H}_v,x),x)
    =&\, x^{m(r-2)}{\varphi^m(\mathcal{H},x)}[x\varphi(\mathcal{G}-u,x)-\sum_{e\in E_{\mathcal{G}}(u)}\varphi\left(\mathcal{G}-V(e),x\right)] \nonumber\\
    & -mx^{(m-1)(r-2)}\varphi(\mathcal{G}-u,x)\varphi(\mathcal{H}-v,x){\varphi^{m-1}(\mathcal{H},x)}.\label{2.4}
      \end{align}
   By replacing
   $x \varphi(\mathcal{G}-u, x)-\sum_{e \in E_{\mathcal{G}}(u)} \varphi(\mathcal{G}-V(e), x)$ by $\varphi(\mathcal{G},x)$ (by Lemma  \ref{su}(b)) in  (\ref{2.4}) and  extracting the common factors
    $x^{(m-1)(r-2)}$ and $\varphi^{m-1}(\mathcal{H},x)$  from the first and the second terms on the right-hand side  of (\ref{2.4}), we get
  \begin{align}
    &\varphi(\mathcal{G}_u\cdot m\mathcal{H}_v,x)\nonumber\\&
    =x^{(m-1)(r-2)}{\varphi^{m-1}({\mathcal{H}},x)}
    \left[x^{r-2}\varphi(\mathcal{G},x)\varphi(\mathcal{H},x)-m\varphi(\mathcal{G}-u,x)\varphi(\mathcal{H}-v,x)\right].\label{2.6}
  \end{align}
  Therefore, by  Lemma \ref{su}(a) and (\ref{2.6}), we have
  \begin{align}
    &\varphi\left[(\mathcal{G}_u\cdot m\mathcal{H}_v,x)\cup (m-1)\mathcal{G},x\right]=
    \varphi(\mathcal{G}_u\cdot m\mathcal{H}_v,x),x){\varphi^{m-1}(\mathcal{G},x)}\nonumber\\
    &= x^{(m-1)(r-2)}{\left[\varphi(\mathcal{G},x)\varphi(\mathcal{H},x)\right]}^{m-1}
    \left[x^{r-2}\varphi(\mathcal{G},x)\varphi(\mathcal{H},x)-m\varphi(\mathcal{G}-u,x)\varphi(\mathcal{H}-v,x)\right].\label{2.8}
  \end{align}
  Similarly, we obtain
  \begin{align}
    &\varphi\left[(\mathcal{H}_v\cdot m \mathcal{G}_u)\cup (m-1)\mathcal{H},x\right]=\varphi(\mathcal{H}_v\cdot m\mathcal{G}_u,x)\varphi^{m-1}( \mathcal{H},x)\nonumber\\
    &=x^{(m-1)(r-2)}{\left[\varphi(\mathcal{G},x)\varphi(\mathcal{H},x)\right]}^{m-1}
    \left[x^{r-2}\varphi(\mathcal{G},x)\varphi(\mathcal{H},x)-m\varphi(\mathcal{G}-u,x)\varphi(\mathcal{H}-v,x)\right].\label{2.8jia}
  \end{align}
  Since the right-hand sides of   (\ref{2.8}) and (\ref{2.8jia}) are the same, it follows from (\ref{2.8}) and (\ref{2.8jia}) that  Theorem \ref{th2} holds.
 ~~$\Box$

  It should be noted that Lemma \ref{th2} generalizes many  known results,
  for example,
  Equation (2a) derived by Cvetkovi\'{c} et al.\ \cite{cvetkovic1975spectral},
  Corollay 2.9 deduced by Shen et al.\ \cite{Shen2005},
   and
  Lemma 3.5(2) obtained by Wang et al.\ \cite{Wang2009}.
  The three results mentioned here are  special cases of Lemma \ref{th2}.

  By Lemmas \ref{large} and \ref{th2},  we obtain Theorem \ref{c4}.

 \begin{theorem} \label{c4}
  Let $m$ be a positive integer. Suppose that $\mathcal{G}$ and $\mathcal{H}$ are two $r$-uniform supertrees with $u\in V(\mathcal{G})$ and $v\in V(\mathcal{H})$, where $r \geq 3$.
  We have (i). $\rho\left(\left(\mathcal{G}_u\cdot m\mathcal{H}_v\right)\cup (m-1)\mathcal{G}\right)
  =\rho\left(\left(\mathcal{H}_v\cdot m\mathcal{G}_u\right)\cup (m-1)\mathcal{H}\right)$
  and $\rho\left(\mathcal{G}_u\cdot m\mathcal{H}_v\right)=\rho\left(\mathcal{H}_v\cdot m\mathcal{G}_u\right)$;
  \textcolor{black}{ (ii). $MG(\left(\mathcal{G}_u\cdot m\mathcal{H}_v\right)\cup \\(m-1)\mathcal{G})
  =MG(\left(\mathcal{H}_v\cdot m\mathcal{G}_u\right)\cup (m-1)\mathcal{H})$}.
  \end{theorem}
%
 \noindent\textbf{Proof}. Obviously,
  $\rho\left(\left(\mathcal{G}_u\cdot m\mathcal{H}_v\right)\cup (m-1)\mathcal{G}\right)=\rho\left(\left(\mathcal{H}_v\cdot m\mathcal{G}_u\right)\cup (m-1)\mathcal{H}\right)$ follows from Lemmas \ref{large} and \ref{th2}.
 Since $\mathcal{G}$ and $\mathcal{H}$ are  respectively the proper subgraphs of
  $\mathcal{G}_u\cdot m\mathcal{H}_v$ and $\mathcal{H}_v\cdot m\mathcal{G}_u$, by Lemma \ref{sub},  we get
 $\rho\left( \mathcal{G}\right)<\rho\left(\mathcal{G}_u\cdot m\mathcal{H}_v\right)$ and $\rho\left( \mathcal{H}\right)<\rho\left(\mathcal{H}_v\cdot m\mathcal{G}_u\right)$.
  It follows from  Lemma \ref{bing} that $\rho\left(\left(\mathcal{G}_u\cdot m\mathcal{H}_v\right)\cup (m-1)\mathcal{G}\right)=\operatorname{max}\{\rho\left(\mathcal{G}_u\cdot m\mathcal{H}_v\right),\rho\left( \mathcal{G}\right)\}=\rho\left(\mathcal{G}_u\cdot m\mathcal{H}_v\right)$. Similarly,  we have
 $\rho\left(\left(\mathcal{H}_v\cdot m\mathcal{G}_u\right)\cup (m-1)\mathcal{H}\right)=\rho\left(\mathcal{H}_v\cdot m\mathcal{G}_u\right)$.
 Thus, we obtain $\rho\left(\mathcal{G}_u\cdot m\mathcal{H}_v\right)=\rho\left(\mathcal{H}_v\cdot m\mathcal{G}_u\right)$. \textcolor{black}{Therefore, we have Theorem \ref{c4}(i).
 By the definition of the matching energy of an $r$-uniform  hypergraph and Lemma \ref{th2}, we obtain Theorem \ref{c4}(ii)}.
 ~~$\Box$

  By Lemmas \ref{tp} and \ref{th2}, we can directly get Theorem 3.2(ii) in Ref. \cite{Wu2018} which was obtained by Wu et al.\ The result  is shown in Theorem \ref{c3s}.
  %

 \begin{theorem}\label{c3s}
  Let $m$ be a positive integer. Suppose that  $G$ and $H$ are two trees with $u\in V(G)$ and $v\in V(H)$.
  Then  $\phi\left[\left(G_u\cdot mH_v\right)\cup (m-1)G,x\right]=\phi\left[\left(H_v\cdot mG_u\right)\cup (m-1)H,x\right]$.
  Namely, $(G_u\cdot mH_v)\cup (m-1)G$ and
  $(H_v\cdot  mG_u)\cup (m-1)H$  are $\bm{M}$-cospectral, where $\bm{M}$ is the adjacency matrix.
   \end{theorem}

 \section{The third pair of  $r$-uniform supertrees with the same spectral radius and matching energy}

In this section, we characterize the third  pair of $r$-uniform supertrees with the same SR \textcolor{black}{and ME,
  and get a  graph which is not determined by its  spectra of its adjacency matrix. The two results are  shown in
 Theorems \ref{c5} and \ref{c6}. To obtain our results, Lemma \ref{th3} is introduced first.}

 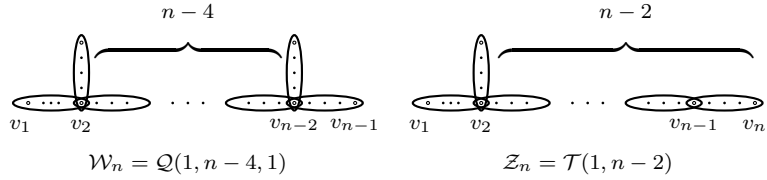
\begin{figure}
\thicklines
\begin{tikzpicture}
\draw[thick,black](2.3,-1.2) ellipse(0.5 and 0.1); 
\draw[black](2.0,-1.2) circle(0.02);
\draw(1.9,-1.5) node{{\scriptsize $v_1$}};
\fill(2.2, -1.2) circle(0.5pt); 
\fill(2.3, -1.2) circle(0.5pt); 
\fill(2.4, -1.2) circle(0.5pt); 
\draw[black](2.7,-1.2) circle(0.02);
\draw(2.7,-1.5) node{{\scriptsize $v_2$}};
\draw[thick,black,rotate around={90:(2.7,-1.2)}](3.1,-1.2) ellipse(0.5 and 0.1);
\fill(2.7, -0.6) circle(0.5pt); 
\fill(2.7, -0.8) circle(0.5pt); 
\fill(2.7, -1.0) circle(0.5pt); 
\draw[black](2.7,-0.4) circle(0.02);
\draw[thick,black](3.1,-1.2) ellipse(0.5 and 0.1); 
\fill(2.9, -1.2) circle(0.5pt); 
\fill(3.1, -1.2) circle(0.5pt); 
\fill(3.3, -1.2) circle(0.5pt); 
\fill(3.9, -1.2) circle(0.5pt); 
\fill(4.1, -1.2) circle(0.5pt); 
\fill(4.3, -1.2) circle(0.5pt); 
\draw[thick,black](5.1,-1.2) ellipse(0.5 and 0.1); 
\fill(4.9, -1.2) circle(0.5pt); 
\fill(5.1, -1.2) circle(0.5pt); 
\fill(5.3, -1.2) circle(0.5pt); 
\draw[black](5.5,-1.2) circle(0.02);
\draw(5.5,-1.5) node{{\scriptsize $v_{n-2}$}};
\draw[thick,black,rotate around={90:(5.5,-1.2)}](5.9,-1.2) ellipse(0.5 and 0.1);
\fill(5.5, -0.6) circle(0.5pt); 
\fill(5.5, -0.8) circle(0.5pt); 
\fill(5.5, -1.0) circle(0.5pt); 
\draw[black](5.5,-0.4) circle(0.02);
\draw(4.1,0) node{{\scriptsize $n-4$}};
\draw(4.1,-0.5) node{$\overbrace{~~~~~~~~~~~~~~~~~~~}$};
\draw[thick,black](5.9,-1.2) ellipse(0.5 and 0.1); 
\fill(5.7, -1.2) circle(0.5pt); 
\fill(5.9, -1.2) circle(0.5pt); 
\fill(6.1, -1.2) circle(0.5pt); 
\draw[black](6.3,-1.2) circle(0.02);
\draw(6.3,-1.5) node{{\scriptsize $v_{n-1}$}};
\draw(4.1,-2) node{{\scriptsize $\mathcal{W}_n=\mathcal{Q}(1,n-4,1)$}};
\end{tikzpicture}
\begin{tikzpicture}
\draw[thick,black](2.3,-1.2) ellipse(0.5 and 0.1); 
\draw[black](2.0,-1.2) circle(0.02);
\draw(1.9,-1.5) node{{\scriptsize $v_1$}};
\fill(2.2, -1.2) circle(0.5pt); 
\fill(2.3, -1.2) circle(0.5pt); 
\fill(2.4, -1.2) circle(0.5pt); 
\draw[black](2.7,-1.2) circle(0.02);
\draw(2.7,-1.5) node{{\scriptsize $v_2$}};
\draw[thick,black,rotate around={90:(2.7,-1.2)}](3.1,-1.2) ellipse(0.5 and 0.1);
\fill(2.7, -0.6) circle(0.5pt); 
\fill(2.7, -0.8) circle(0.5pt); 
\fill(2.7, -1.0) circle(0.5pt); 
\draw[black](2.7,-0.4) circle(0.02);
\draw[thick,black](3.1,-1.2) ellipse(0.5 and 0.1); 
\fill(2.9, -1.2) circle(0.5pt); 
\fill(3.1, -1.2) circle(0.5pt); 
\fill(3.3, -1.2) circle(0.5pt); 
\fill(3.9, -1.2) circle(0.5pt); 
\fill(4.1, -1.2) circle(0.5pt); 
\fill(4.3, -1.2) circle(0.5pt); 
\draw[thick,black](5.1,-1.2) ellipse(0.5 and 0.1); 
\fill(4.9, -1.2) circle(0.5pt); 
\fill(5.1, -1.2) circle(0.5pt); 
\fill(5.3, -1.2) circle(0.5pt); 
\draw[black](5.5,-1.2) circle(0.02);
\draw(5.5,-1.5) node{{\scriptsize $v_{n-1}$}};
\draw(4.6,0) node{{\scriptsize $n-2$}};
\draw(4.6,-0.5) node{$\overbrace{~~~~~~~~~~~~~~~~~~~~~~~~~~}$};
\draw[thick,black](5.9,-1.2) ellipse(0.5 and 0.1); 
\fill(5.7, -1.2) circle(0.5pt); 
\fill(5.9, -1.2) circle(0.5pt); 
\fill(6.1, -1.2) circle(0.5pt); 
\draw[black](6.3,-1.2) circle(0.02);
\draw(6.3,-1.5) node{{\scriptsize $v_{n}$}};
\draw(4.1,-2) node{{\scriptsize $\mathcal{Z}_n=\mathcal{T}(1,n-2)$}};
\end{tikzpicture}
\caption{\label{c3} $\mathcal{W}_n=\mathcal{Q}(1,n-4,1)$ and $\mathcal{Z}_n=\mathcal{T}(1,n-2)$}
\end{figure}

 Let $\mathcal{W}_n=\mathcal{Q}(1,n-4,1)$ with $n\geq 5$ and $\mathcal{Z}_n=\mathcal{T}(1,n-2)$ with $n\geq 2$.
  $\mathcal{W}_n$ and $\mathcal{Z}_n$ are shown in Fig.\ \ref{c3}.
  The loose path in $\mathcal{W}_{n}$ is denoted  by $\mathcal{P}_{n-2}^r=v_1e_1v_2e_2v_3\cdots v_{n-2}e_{n-2}v_{n-1}$, where $n\geq 5$.

\begin{lemma}\label{th3}
  We have 
  $\varphi(\mathcal{P}_{m-5}^r\cup\mathcal{W}_{n-1},x)=\varphi(\mathcal{P}_{n-5}^r\cup\mathcal{W}_{m-1},x)$, where  $m, n\geq6$ and $r\geq 2$.
  \end{lemma}

   \noindent\textbf{Proof}.
   (1). The proof of Lemma  \ref{th3} when $m=6$ and $n\geq 6$.

   When $m=6$ and $n\geq 6$, we prove
   \begin{align}\label{www}
   \varphi(\mathcal{P}_{1}^r\cup\mathcal{W}_{n-1},x)=\varphi(\mathcal{P}_{n-5}^r\cup\mathcal{W}_5,x)
   \end{align}
   by induction on $n$.

  (i). When $n=6$, obviously, $\varphi(\mathcal{P}_1^r\cup\mathcal{W}_5,x)=\varphi(\mathcal{P}_1^r\cup\mathcal{W}_5,x)$.

   When  $n=7$, by Lemma \ref{redf}, we get
    $\varphi(\mathcal{P}_1^r,x)=x^r-1$,
    $\varphi(\mathcal{W}_6,x)=\varphi(\mathcal{Q}(1,2,1),x)=x^{6r-5}-6x^{5r-5}+8x^{4r-5}$,
    $\varphi(\mathcal{P}_2^r,x)=x^{2r-1}-2x^{r-1}$, and
    $\varphi(\mathcal{W}_5,x)=\varphi(\mathcal{Q}(1,1,1),x)=x^{5r-4}-5x^{4r-4}+4x^{3r-4}$.
   Thus, we get
  \begin{align*}
  \varphi(\mathcal{P}_1^r\cup\mathcal{W}_6,x)&=\varphi(\mathcal{P}_2^r\cup\mathcal{W}_5,x)=x^{7r-5}-7x^{6r-5}+14x^{5r-5}-8x^{4r-5}.
  \end{align*}

 (ii).  When $n=k$ with $k\geq 7$, we suppose that (\ref{www}) hold. Namely, $\varphi(\mathcal{P}_{1}^r\cup\mathcal{W}_{k-1},x)=\varphi(\mathcal{P}_{k-5}^r\cup\mathcal{W}_5,x)$.

  (iii).  When $n=k+1$ with $k\geq 7$, we  prove that (\ref{www}) hold.

  We have
   \begin{align}
 \varphi(\mathcal{W}_{n},x) &=x\varphi(\mathcal{W}_{n}-v_{n-1},x)-\sum_{e\in E_{\mathcal{W}_{n}}(v_{n-1})}\varphi(\mathcal{W}_{n}-V(e),x)\label{3.1}\\
  &=x\varphi(\mathcal{Z}_{n-1}\cup N_{r-2},x)-\varphi(\mathcal{Z}_{n-3}\cup N_{r-1}\cup N_{r-2},x)\label{3.2}\\
  &=x^{r-1}[\varphi(\mathcal{Z}_{n-1},x)-x^{r-2}\varphi(\mathcal{Z}_{n-3},x)],\label{3.44}
\end{align}
 where \eqref{3.1} follows from Lemma \ref{su}(b),
 \eqref{3.2} holds since $\mathcal{W}_{n}-v_{n-1}\cong \mathcal{Z}_{n-1}\cup N_{r-2}$ and  $E_{\mathcal{W}_{n}}(v_{n-1})=\{e_{n-2}\}$,
 and \eqref{3.44} is derived from Lemma \ref{su}(a).
 Similarly, we get
 \begin{align}\label{3.9}
  \varphi(\mathcal{W}_{n-1},x)
  &=x^{r-1}[\varphi(\mathcal{Z}_{n-2},x)-x^{r-2}\varphi(\mathcal{Z}_{n-4},x)],
 \end{align}
 \begin{align}\label{3.10}
  \varphi(\mathcal{W}_{n-2},x)
  &=x^{r-1}[\varphi(\mathcal{Z}_{n-3},x)-x^{r-2}\varphi(\mathcal{Z}_{n-5},x)].
 \end{align}

  Let the loose path of $\mathcal{Z}_{n-1}$ be $\mathcal{P}_{n-2}^r=v'_1e'_1v'_2e'_2v'_3\cdots v'_{n-2}e'_{n-2}v'_{n-1}$.
  We obtain
 \begin{align}
  \varphi(\mathcal{Z}_{n-1},x)&=x\varphi(\mathcal{Z}_{n-1}-v'_{n-1},x)-\sum_{e\in E_{\mathcal{Z}_{n-1}}(v'_{n-1})}\varphi(\mathcal{Z}_{n-1}-V(e),x)\nonumber\\
  &=x\varphi(\mathcal{Z}_{n-2}\cup N_{r-2},x)-\varphi(\mathcal{Z}_{n-3}\cup N_{r-2},x)\nonumber\\
  &=x^{r-2}\left[x\varphi(\mathcal{Z}_{n-2},x)-\varphi(\mathcal{Z}_{n-3},x)\right]\label{3.4}
 \end{align}
 Similarly, we have
 \begin{align}\label{33.5}
  \varphi(\mathcal{Z}_{n-3},x)&=x^{r-2}\left[x\varphi(\mathcal{Z}_{n-4},x)-\varphi(\mathcal{Z}_{n-5},x)\right].
 \end{align}
 Thus, by substituting  \eqref{3.4} and \eqref{33.5}  into \eqref{3.44}, we obtain
 \begin{align}
   &\varphi(\mathcal{W}_{n},x)=x^{r-1}[\varphi(\mathcal{Z}_{n-1},x)-x^{r-2}\varphi(\mathcal{Z}_{n-3},x)]\nonumber\\
   &=x^{r-1}\!\bigg[x^{r-2}\left[x\varphi(\mathcal{Z}_{n-2},x)-\varphi(\mathcal{Z}_{n-3},x)\right]
   -x^{r-2}\Big[x^{r-2}\left[x\varphi(\mathcal{Z}_{n-4},x)-\varphi(\mathcal{Z}_{n-5},x)\right]\!\Big]\!\bigg]\label{3.11}\\
   &=x^{r-2}\Big[x^r[\varphi(\mathcal{Z}_{n-2},x)-x^{r-2}
   \varphi(\mathcal{Z}_{n-4},x)]-x^{r-1}[\varphi(\mathcal{Z}_{n-3},x)-x^{r-2}\varphi(\mathcal{Z}_{n-5},x)]\Big]\nonumber\\
   &=x^{r-2}[x\varphi(\mathcal{W}_{n-1},x)-\varphi(\mathcal{W}_{n-2},x)].  \label{3.1111}
 \end{align}
 It is noted that 
 \eqref{3.1111} is derived from \eqref{3.9} and \eqref{3.10}.
  Let $k\geq 7$. Thus, we obtain
 \begin{align}
   \varphi(\mathcal{P}_1^r\cup\mathcal{W}_{k},x)
   &=\varphi(\mathcal{P}_1^r,x)\varphi(\mathcal{W}_{k},x)\label{3.12}\\
   &=\varphi(\mathcal{P}_1^r,x)x^{r-2}\left[x\varphi(\mathcal{W}_{k-1},x)-\varphi(\mathcal{W}_{k-2},x)\right]\label{3.13}\\
   &=x^{r-2}\left[x\varphi(\mathcal{P}_1^r\cup\mathcal{W}_{k-1},x)-\varphi(\mathcal{P}_1^r\cup\mathcal{W}_{k-2},x)\right]\label{3.14}\\
   &=x^{r-2}\left[x\varphi(\mathcal{P}_{k-5}^r\cup\mathcal{W}_5,x)-\varphi(\mathcal{P}_{k-6}^r\cup\mathcal{W}_5,x)\right]\label{3.15}\\
   &=x^{r-2}\varphi(\mathcal{W}_5,x)\left[x\varphi(\mathcal{P}_{k-5}^r,x)-\varphi(\mathcal{P}_{k-6}^r,x)\right].\label{3.16}
 \end{align}
 It is noted that  \eqref{3.12} follows from Lemma \ref{su}(a),
 \eqref{3.13} is obtained by substituting  \eqref{3.1111}  into \eqref{3.12},
 \eqref{3.14} is derived from Lemma \ref{su}(a),
 \eqref{3.15} follows from  $\varphi(\mathcal{P}_1^r\cup\mathcal{W}_{k-2},x)=\varphi(\mathcal{P}_{k-6}^r\cup\mathcal{W}_5,x)$ and
 $\varphi(\mathcal{P}_1^r\cup\mathcal{W}_{k-1},x)=\varphi(\mathcal{P}_{k-5}^r\cup\mathcal{W}_5,x)$
 (by the inductive hypothesis), and
 \eqref{3.16} is deduced from Lemma \ref{su}(a).
 Furthermore, we get
\begin{align}\label{3.17}
   \varphi(\mathcal{P}_{k-4}^r,x)
   &=x^{r-2}\left[x\varphi(\mathcal{P}_{k-5}^r,x)-\varphi(\mathcal{P}_{k-6}^r,x)\right].
 \end{align}
 By substituting  \eqref{3.17} into \eqref{3.16}, we have
 \begin{align}\label{3.18}
   \varphi(\mathcal{P}_1^r\cup\mathcal{W}_{k},x)&=\varphi(\mathcal{W}_5,x)\varphi(\mathcal{P}_{k-4}^r,x)=\varphi( \mathcal{P}_{k-4}^r \cup \mathcal{W}_5,x).
 \end{align}
 Therefore, by the method of inductive  hypothesis, when $m=6$ and $n\geq 6$, we obtain \eqref{www}. 

  (2). The proof of Lemma  \ref{th3} when $m\geq 7$ and $n\geq 6$.

  Let $m\geq 7$ and $n\geq 6$. By (\ref{www}) and Lemma \ref{su}(a), we have
  $\varphi(\mathcal{P}_1^r\cup\mathcal{W}_{n-1}\cup \mathcal{P}_{m-5}^r,x)=\varphi(\mathcal{P}_{n-5}^r\cup\mathcal{W}_5\cup \mathcal{P}_{m-5}^r,x)=\varphi(\mathcal{P}_{n-5}^r\cup\mathcal{W}_{m-1}\cup \mathcal{P}_1^r,x)$.
  Therefore, by  Lemma \ref{su}(a), we get
   $ \varphi(\mathcal{P}_{m-5}^r\cup\mathcal{W}_{n-1},x)=\varphi(\mathcal{P}_{n-5}^r\cup\mathcal{W}_{m-1},x)$.

   By combining the proofs of (1) and (2), we get Lemma \ref{th3}.
    ~~$\Box$

    By Lemmas \ref{large} and \ref{th3}, in Theorem \ref{c5},  we obtain the third  pair  of $r$-uniform supertrees with the same SR \textcolor{black}{and ME}, where $r\geq 3$.

  \begin{theorem}\label{c5}
  Let $m, n\geq6$ and $r\geq 3$. We have
   $\rho(\mathcal{P}_{m-5}^r\cup\mathcal{W}_{n-1})=\rho(\mathcal{P}_{n-5}^r\cup\mathcal{W}_{m-1})$ and
   \textcolor{black}{$ME(\mathcal{P}_{m-5}^r\cup\mathcal{W}_{n-1})=ME(\mathcal{P}_{n-5}^r\cup\mathcal{W}_{m-1})$}.
     \end{theorem}


  In the following, all mentioned results are related with the spectra of adjacency matrix.
 Let  $W_n=Q(1,n-4, 1)$ with $n\geq 5$ and $Z_n=T(1,n-2)$ with $n\geq 2$.
 Shen et al.\ \cite{Shen2005} deduced that  $P_{n-1}\cup Z_{n+1}$ ($n\geq 1$) is not determined by their  spectra while $Z_{n+1}$ ($n\geq 1$) and $Z_{n_1+1}\cup\cdots\cup Z_{n_k+1}$ ($n_1,n_2,\ldots,n_k \geq 2$) are determined by their  spectra.
  Wang et al.\ \cite{Wang2009} obtained that $P_0\cup{P_{n-1}}$ is determined by their  spectra if and only if $n=2k$ with $k\geq 1$.
 Cvetkovi\'{c} and Jovanovi{\'c} \cite{Cvetkovic2017} derived that $Z_{n+1}\cup{P_0}$ ($n\geq9$) is determined by their   spectra.
 Inspired by all the above-mentioned results,
 in Theorem \ref{c6}, we get that $W_{n-1}\cup P_{m-5}$ ($n,m\geq6$ and $n\neq m$) is not determined by its spectra.


  \begin{theorem}\label{c6}
  Let $m, n\geq6$. We have
  $\phi(P_{m-5}\cup W_{n-1},x)=\phi(P_{n-5}\cup W_{m-1},x)$.
 Namely,  $P_{m-5}\cup W_{n-1}$ and $P_{n-5}\cup W_{m-1}$ are $\bm{M}$-cospectral and \textcolor{black}{$ME(P_{m-5}\cup W_{n-1})=ME(P_{n-5}\cup W_{m-1})$}, where  $\bm{M}$ is the adjacency matrix.
    \end{theorem}

  \noindent\textbf{Proof}. Let $m, n\geq6$.
      By  Lemma  \ref{th3}, when $r=2$,  we get  $\varphi(P_{m-5}\cup W_{n-1},x)=\varphi(P_{n-5}\cup W_{m-1},x)$.
         Furthermore,  by  Lemma \ref{tp} and the definition of the matching energy of a graph, we get Theorem \ref{c6}.
           ~~$\Box$

%

%
\begin{acknowledgments}
\noindent\textbf{Acknowledgments}\\
 The work was supported by the Natural Science Foundation of Shanghai under the grant number 21ZR1423500 and the National Natural Science Foundation of China under the grant number 11871040.
\end{acknowledgments}

\bibliographystyle{elsevier_citation_order}  
\bibliography{zjx20211223}

\end{document}